\documentclass[11pt]{amsart}
\usepackage{amssymb,amsmath,amsthm,amsfonts,eucal,amscd}

\newcommand{\ol}{\overline}

\theoremstyle{plain}
\newtheorem{bigthm}{\sc Theorem}
\newtheorem{Theorem}{\sc Theorem}[section]
\newtheorem{Lemma}[Theorem]{\sc Lemma}
\newtheorem{Proposition}[Theorem]{\sc Proposition}
\newtheorem{Corollary}[Theorem]{\sc Corollary}

\newtheorem*{Acknowledgement}{\sc Acknowledgement}

\theoremstyle{definition}
\newtheorem{Definition}{\sc Definition}[section]

\theoremstyle{remark}

\newtheorem{Remark}{\sc Remark}[section]

\numberwithin{equation}{section}

\newcommand{\bm}[1]{\mbox{$\mathbf{#1}$}}%

\def\ol{\overline}%
\def\op{{operator~}}%
\def\ops{{operators~}}%
\def\hs{{Hilbert space~}}%

\newcommand{\rk}{reproducing kernel }%

\newcommand{\qand}{\quad\text{and}\quad}
\newcommand{\qfor}{\quad\text{for}\quad}

\def\C{\mbox{${\mathbb C}$}}

\def\B{\mbox{${\mathbb B}$}}
\def\D{\mbox{${\mathbb D}$}}

\def\N{\mbox{${\mathbb N}$}}

\def\F{\mbox{${\mathbb F}$}}
\def\1{\mbox{${\mathbb 1}$}}
\def\I{\mbox{${\mathbb I}$}}

\newcommand{\scZ}{\mathcal{Z}}

\newcommand{\scR}{\mathcal{R}}
\newcommand{\scJ}{\mathcal{J}}

\newcommand{\scH}{\mathcal{H}}
\newcommand{\scK}{\mathcal{K}}
\newcommand{\scM}{\mathcal{M}}

\newcommand{\scQ}{\mathcal{Q}}
\newcommand{\scA}{\mathcal{A}}
\newcommand{\scS}{\mathcal{S}}

\newcommand{\scD}{\mathcal{D}}
\newcommand{\scE}{\mathcal{E}}
\newcommand{\scT}{\mathcal{T}}

\newcommand{\alg}[1]{\mbox{${\mathcal A}(#1)$}}%

\newcommand{\exactseq}[3]{\mbox{$0 \longleftarrow  {#1} \twoheadleftarrow
{#2} \hookleftarrow {#3} \longleftarrow 0$}}%

\def\x{\mbox{{$\mathbf{x}$}}}%
\def\h{\mbox{{$\mathbf{h}$}}}%
\def\z{\mbox{{$z^\prime$}}}%
\def\w{\mbox{{$\mathbf{w}$}}}%
\def\0{\mbox{{\Large$\mathbf{0}$}}}%

\newcommand{\vect}[2]{\begin{pmatrix} #1 \\ #2 \end{pmatrix}}%
\newcommand{\kf}[1]{{K(\cdot,#1)}}%
\newcommand{\ktf}[2]{{K(#1,#2)}}%

\newcommand{\deru}{\mbox{$\partial_1 \bar{\partial}_1$}}%

\newcommand{\der}{\mbox{$\frac{\partial^2}{\partial w\partial\bar{w}}$}}%

\newcommand{\inner}[2]{\langle #1,#2 \rangle }%
\newcommand{\norm}[1]{\mbox{$\left \| #1\right \|$}}%

\title{Equivalence of quotient Hilbert modules -- II}
\author{Ronald G. Douglas}
\address{Texas A\&M University\\College Station\\ Texas 77843-3368 }
\email[Ronald G. Douglas]{rdouglas@math.tamu.edu}
\author{Gadadhar Misra}
\address{Indian Statistical Institute\\
R. V. College Post\\
Bangalore 560 059 } \email[Gadadhar Misra]{gm@isibang.ac.in}
\subjclass[1991]{46E22, 32Axx, 32Qxx, 47A20, 47A65, 47B32 and 55R65} 
\keywords{Hilbert modules, Complex geometry, Jet bundles, Curvature, Homogeneous operators}

\thanks{The research of both the authors was supported in part by a
grant from the DST - NSF  Science and Technology Cooperation Programme.}

\begin{document}

\begin{abstract}

  For any open, connected and bounded set $\Omega\subseteq \C^m$, let
  $\scA$ be a natural function algebra consisting of functions
  holomorphic on $\Omega$.  Let $\scM$ be a Hilbert module over the
  algebra $\scA$ and $\scM_0\subseteq \scM$ be the submodule of
  functions vanishing to order $k$ on a hypersurface $\scZ \subseteq
  \Omega$.  Recently the authors have obtained an explicit
  complete set of unitary invariants for the quotient module
  $\scQ=\scM \ominus \scM_0$ in the case of $k=2$.  In this paper, we
  relate these invariants to familiar notions from complex geometry.
  We also find a complete set of unitary invariants for the general
  case.  We discuss many concrete examples in this setting.  As an
  application of our equivalence results, we characterise certain 
  homogeneous Hilbert modules over the bi-disc algebra.

\end{abstract}

\maketitle



\section{Introduction}\label{}

One source of fascination in the study of operator theory is the wide
variety of connections made with other branches of mathematics.
Techniques from algebra, topology, geometry and analysis are used to
understand bounded linear operators on Hilbert space.  And, in many
instances, the behavior and properties of the operators can be used to
illustrate critical features and aspects of the other fields.  This is
particularly true in the case of multivariate operator theory, that
is, when several operators or an algebra of operators is studied.
Here the setting and results from these other areas can be quite
sophisticated and the techniques used to understand multivariate
operator theory often require additional development.  Such is the
focus of this paper.

Although the spectral theorem is a key tool in the study of
self-adjoint and normal operators, there are large and important
classes of naturally occurring operators to which this theory doesn't
apply. Examples illustrating such phenomena can be obtained by
considering multiplication operators on spaces of holomorphic
functions on some domain in $\C^m$.  For domains in $\C$, one is in the
realm of single operator theory while it is multivariate operator
theory for $m > 1$.  If one considers the unit ball $\B^m$ in $\C^m$ and the
Bergman space $A^2(\B^m)$ for it, one obtains a module over the polynomial
algebra $\C[z]$, where $z = (z_1,\ldots, z_m)$.

Techniques from complex geometry were shown in \cite{CD},
\cite{Bolyai}, and \cite{CD85} to be useful in studying such Hilbert
modules.  Closed submodules related to polynomial ideals were shown to
reflect properties of the ideals and results in a rigidity phenomenon
for such submodules \cite{D-P-S-Y}.  In \cite{rgdgm}, it was shown
that the study of quotient modules, determined by polynomial ideals,
could also be reduced to the earlier work involving complex geometry
if the ideal is principal and prime.  As might be expected, the
non-prime case is more complicated (cf. \cite{DMV}) and some real
technical difficulties arise in the complex geometry needed to handle
its study.  Overcoming these problems by developing new results and
techniques in complex geometry is the main goal of this paper.

In a basic construction, Hilbert modules, such as $A^2(\B^m)$, can be
shown to yield a hermitian holomorphic vector bundle over the domain
and this bundle characterizes the module up to unitary equivalence.
Moreover, the geometric invariants of the bundle, including the
curvature, can be obtained from the module action.
For quotient modules by multiplicity-free principal ideals, a bundle
still exists but over the intersection of the domain with the zero
variety of the ideal. One can also exhibit a kernel
function that characterizess the quotient module. 
In \cite[Theorem 1.4]{rgdgm},
the fundamental class of the hypersurface $\scZ$ was expressed using
the curvatures of the pair of modules $\scM$ and the submodule
$\scM_0$ and the localization of the inclusion map $\scM_0
\hookrightarrow \scM$. Here we give a complete set of invariants for
the equivalence of the quotient modules. In \cite{DMV} it was shown that for quotient
modules obtained using submodules of higher multiplicity, there is still a
higher rank hermitian holomorphic bundle but it's bundle structure is
not enough to characterize the quotient module.  One must also involve
the flag structure in the bundle defined by the module action which
now involves nilpotents.  But even that is not enough.  In particular,
we must consider the nilpotent structure defined by
the module action itself.  Classifying such
objects requires introducing new ideas and techniques extending older
ones from complex geometry which involve jet bundles and moving
frames. Before describing our main results, we neeed to introduce 
some terminology.

For any bounded open connected subset $\Omega$ of
$\C^m$, let $\alg{\Omega}$ be the completion, with respect to the
supremum norm on the closure $\ol{\Omega}$ of the domain $\Omega$, of
functions holomorphic in a neighbourhood of $\ol{\Omega}$.  The
Hilbert space $\scM$ is said to be a {\em Hilbert module\/} over
$\scA(\Omega)$ if $\scM$ is a module over $\scA(\Omega)$ with module
map $\scA(\Omega)\times \scM\to \scM$ defined by pointwise
multiplication such that
$$
\|f \cdot h\|_{\scM} \le C \|f\|_{\scA(\Omega)} \|h\|_{\scM}
\qfor f \in \scA(\Omega) \qand h \in \scM,
$$
for some positive constant $C$ independent of $f$ and $h$.
It is said to be {\em contractive \/} if we also have
$C \leq 1$.


We work in a class of locally free Hilbert modules
called quasifree which is defined in Section 2. 
Now fix a hypersurface $\mathcal Z \subseteq \Omega$ 
and let $\scM_0 \subseteq \scM$ be the submodule of the quasi-free
Hilbert module $\scM$ of rank $1$ consisting of those functions in
$\scM$ which vanish to order $k$ on the hypersurface $\scZ$. The
quotient $\scQ = \scM\ominus \scM_0$ is a Hilbert module over
$\alg{\Omega}$, where the module action is naturally defined as $f
\cdot (h + \scM_0) = f\cdot h + \scM_0$ (cf. \cite[Definition
2.2]{rgdvip}.  In other words,
\begin{equation} \label{quotient}
\exactseq{\scQ}{\scM}{\scM_0},
\end{equation}
is an exact sequence of Hilbert modules over $\alg{\Omega}$, where
$\twoheadleftarrow$ is the quotient map and $\hookleftarrow$ is the inclusion map.
It is then possible to obtain geometric invariants for the quotient module
$\scQ$ using the module map $\scM \twoheadleftarrow \scM_0$  
(cf. \cite[Theorem 1.4]{rgdgm}).

For any fixed but arbitrary 
$u\in \Omega$, we may pick a small enough open neighborhood $U
\subseteq \Omega$ of $u$ such that $U \cap \mathcal Z$ admits a
defining function, say $\varphi$, with gradient of $\varphi$ not zero
in the normal direction to $U \cap \mathcal Z$.  Since $U$ is open in $\Omega$, it follows
that the module $\scM$ and the submodule $\scM_0 \subseteq \scM$ are
isomorphic to the module $\scM_{|{\rm res~}U}$ and the submodule
$(\scM_0)_{|{\rm res~}U}\subseteq \scM_{|{\rm res~}U}$ of functions in
$\scM_{|{\rm res~}U}$ which vanish on $U \cap \scZ$,
respectively. Consequently, if we choose to work with the latter pair
of modules, then the corresponding quotient module $\scM_{|{\rm
res~}U} \ominus (\scM_0)_{|{\rm res~}U}$ is isomorphic to the quotient
module $\scM\ominus \scM_0$.  {\em Therefore we may cut down, if necessary,  
the domain $\Omega$ to a suitable small open subset 
$U\subseteq \Omega$ and work with the smaller open set $U$ and the hypersurface
$U \cap \scZ \subseteq U$ without loss of generality.}

The submodule ${\scM}_0$ in \cite{DMV} is taken to be the
(maximal) set of functions which vanish to some given order $k$ on
the hypersurface ${\scZ}$. As in the case of multiplicity one, two
descriptions are provided for the quotient module. A matrix --
valued kernel function must now be used and, in the vector bundle
picture, we have a rank $k$ hermitian, anti-holomorphic vector
bundle over $\mathcal Z^*$.  Some invariants for the quotient
module (though not a complete set) are described in \cite{DMV}.

Finally, in the paper \cite{eqhm}, a complete set of unitary
invariants is obtained for the particular case where $\scM_0$
consists of functions vanishing to order $2$ on the hypersurface.
While two of the three invariants obtained there consist of coefficients of
the curvature form for the jet bundle for $E$, the third, which we
called the ``angle'', seemed to be not so familiar.  In this note,
we show that the angle invariant can be replaced by the second
fundamental form corresponding to the inclusion of the bundle $E$
in the jet bundle $J^{(2)}E$.

Our main goal, however, is to obtain invariants that are
complete, computable and natural in complex geometry -- for general
$k$, not just in the case $k=2$.  We'll state our results in the
following section after we introduce the necessary notation. 

We provide a number of applications of our results to the context 
of homogeneous operators. Moreover, we discuss various aspects of 
global versus local differences in connection with the jet bundle 
construction.  Finally, we describe some relations between the new 
invariants we introduce and several other topics relating to the 
moving frames of Cartan and integrability conditions in 
Chern-Moser theory.

In a recent paper \cite{quasi}, we pointed out that much of
what we proved there was valid for modules over algebras of
holomorphic functions that are not complete; for example, C[z] or
the functions holomorphic on the closure of the domain.  This
continues to be the case in this paper as well. However, we will 
continue to state our principal results for Hilbert modules over 
the complete function algebra $\mathcal A(\Omega)$ even though 
a more general statement would be possible.

\section{Main Results}

\subsection{} \label{wells} We recall some basic notions from 
complex geometry following Wells \cite[Chapter III]{RW}.  
If $E$ is a  hermitian holomorphic vector bundle, then there is a canonical 
connection $D$ on the bundle $E$ which is compatible with both the holomorphic 
and hermitian structures. The curvature $\mathcal K$ of the bundle 
$E$ is then simply defined to be $D\circ D$. Let us provide some more details.

Let $E$ be a hermitian holomorphic vector bundle of rank $k$ over
a complex manifold $M$ and $C^\infty_p(M,E)$ be the space of
smooth $p$ - forms on $M$ whose coefficients are smooth sections
of $E$. A connection on the bundle $E$ is a differential operator
$D: C_p^\infty(M,E) \to C_{p+1}^\infty(M,E)$ of order $1$
satisfying
\begin{equation} \label{compmetric}
D(f\wedge s) = df \wedge s + (-1)^p f \wedge Ds
\end{equation}
for any $f\in C_p^\infty(M,\C)$ and $s\in C_p^\infty(M,E)$, where
$df$ stands for the usual exterior derivative of $f$. 

Assume that
$\theta:E_{|U} \to U\times\C^k$ is a trivialization of $E$ over
some open subset $U$ of $M$.  Let $(s_1,\ldots, s_k)$ be the
corresponding frame of $E_{|U}$.  Then any $s\in C^\infty_p(M,E)$
can be written uniquely as
$$s=\sum_j \sigma_j\otimes s_j,~~\sigma_j\in C_p^\infty(U,\C),~~1\leq j\leq k.$$
Using the hermitian structure $h$ of $E$, we can define a natural
sesqui-linear map
\begin{eqnarray*}
\lefteqn{C_p^\infty(M,E) \times C_q^\infty(M,E) \to C_{p+q}^\infty(M,\C)}\\
& & (s_1,s_2) \mapsto \{s_1,s_2\}
\end{eqnarray*}
combining the wedge product of forms with the hermitian metric on
$E$. If $s = \sum_j \sigma_j \otimes s_j$ and $\tilde{s} = \sum_j
\tilde{\sigma}_j \otimes s_j$, then
$$\{s,\tilde{s}\} = \sum_{j,\ell} \sigma_j \wedge \bar{\tilde{\sigma}}_\ell h(s_j,s_\ell).$$

The curvature tensor $\scK$ associated with the canonical
connection $D$ is in $C^\infty_{1,1}(M, \mbox{herm}(E,E))$.
(Here $C^\infty_{p,q}$ represents the space of forms of degree 
$p$ in the holomorphic differentials and $q$ in the anti-holomorphic ones.)
Moreover, if $h$ is a local representation of the metric in some
open set, then the curvature $\scK = \bar{\partial}({h}^{-1} \partial{h})$
\cite[page 82]{RW}. 

\subsection{} To study quotient Hilbert modules determined by a hypesurface, 
we must consider the behavior of the holomorphic tangent bundle to $\Omega$ 
relative to an analytic hypersurface.  The holomorphic tangent bundle
$T\Omega_{|{\rm res~ } \scZ}$ naturally splits as $T\scZ \dot{+}
N\scZ$, where $N\scZ$ is the normal bundle. It can be identified with
the quotient $T\Omega_{|{\rm res~ } \scZ}/T\scZ$. The co-normal bundle
$N^*\scZ$ is the dual of $N\scZ$; it is the sub-bundle of
$T^*\Omega_{|{\rm res~ } \scZ}$ consisting of cotangent vectors that
vanish on $T\scZ \subseteq T\Omega_{|{\rm res ~} \scZ}$.  Indeed,
there is an easy formula for the co-normal bundle of a smooth
hypersurface which we describe now following \cite[page 145]{G-H}.

Suppose $\scZ$ is given by local defining functions
$\varphi_z$ on $U_z\subseteq \Omega$, $z\in \Omega$,  as in
Definition \ref{hypersurface}. The line bundle $[\scZ]$ defined on
$\Omega$ is then given by transition functions $\{\psi_{z w}=
\tfrac{\varphi_z}{\varphi_w}:z,w \in \Omega \}$ on $U_z \cap U_w$.
By definition, $\varphi_z \equiv 0$ on $U_z \cap \scZ$. It follows
that the differential $d\varphi_z$ is a section of the co-normal
bundle $N^*\scZ$. Besides, $d\varphi_z$ is holomorphic and nonzero
everywhere. On $U_z\cap U_w\cap \scZ$, we have $d\varphi_z =
\psi_{zw} d\varphi_w$, that is, $d\varphi_z$ defines a nonzero
global section of the bundle $N^*\scZ\otimes[\scZ]$. Thus
$N^*\scZ\otimes[\scZ]$ is the trivial bundle which gives the
formula $N^*\scZ = [-\scZ]_{|{\rm res ~} \scZ}$, where $[-\scZ]$
is the inverse of the line bundle $[\scZ]$. This is the Adjunction
Formula I \cite[page 146]{G-H}.

In the following calculation, we assume that $\scZ=\{z_1 = 0\}$.  We
show in subsection \ref{z_1} that there is no loss of generality in
doing so. Let $P_1: T^*\Omega_{|{\rm res~} \mathcal Z} \to N^*\scZ$ be
the bundle map which is the projection onto $N^*\scZ$ and
$P_2=(1-P_1):T^*\Omega_{|{\rm res~} \mathcal Z} \to T^*\scZ$ be the
bundle map which is the projection onto $T^*\scZ$. Now, we have a
splitting of the $(1,1)$ forms as follows:
$$
\wedge^{(1,1)}T^* \Omega_{|{\rm res ~}\scZ} = \sum_{i,j=1}^2 P_i\big
(\wedge^{(1,0)}T^*\Omega _{|{\rm res ~}\scZ}\big ) \wedge P_j \big
(\wedge^{(0,1)}T^*\Omega_{|{\rm res ~}\scZ} \big ).
$$
Accordingly, we have the component of the curvature along the
transverse direction to $\scZ$ which we denote by $\scK_{\rm
trans}$. Clearly, $\scK_{\rm trans}= (P_1 \otimes I) \scK_{|{\rm
res~ }\scZ}$. Similarly, let the component of the curvature along
tangential directions to $\scZ$ be $\scK_{\rm tan}$.  Again, $
\scK_{\rm tan} = (P_2 \otimes I)\scK_{|{\rm res~} \scZ}$. (Here
$I$ is the identity map on the vector space $\mbox{herm}(E,E)$.)
In local coordinates, the curvature of $E$ at $ w\in \Omega^*$ is
given by
\begin{eqnarray} \label{curvature}
\scK(w) &=& - \sum_{i,j=1}^m \bar{\partial}_i\Big(K(w,w)^{-1}
{\partial}_j K(w,w) \Big ) d\bar{w}_i\wedge dw_j \nonumber\\
&=& \scK_{11}(w) d\bar{w}_1\wedge dw_1 + \sum_{j=2}^m \scK_{1j}(w) d\bar{w}_1\wedge dw_j
\nonumber \\
&& \phantom{~~~~~} + \sum_{i=2}^m\scK_{i1}(w) d\bar{w}_i\wedge dw_1 +
\sum_{i,j=2}^m\scK_{ij}(w) d\bar{w}_i\wedge dw_j \nonumber\\
&=& \begin{pmatrix} \scK_{\rm tan}(w) & \scS(w) \\
-\overline{\scS(w)} & \scK_{\rm trans}(w) \end{pmatrix}
\vect{{d^\prime}^{\rm t}\bar{w}} {d\bar{w}_1} \wedge \vect{{d^\prime}^{\rm t} w}{d w_1},
\end{eqnarray}
where $\partial_i= \frac{\partial}{\partial w_i}$,
$\bar{\partial}_j= \frac{\partial}{\partial \bar{w}_j}$ and ${d^\prime}^{\rm t} w$ 
denotes the transpose of  $(dw_m, \hdots , dw_2)$. 
Also, we let $\scS(w)$, which appears in the $(1,2)$
position of the decomposition for the curvature, denote the
$(1,1)$ form $\sum_{j=2}^m \scK_{1j}(w) d\bar{w}_1\wedge dw_j$.

We will study quotient modules for a special class of Hilbert modules
which includes the Hardy and Bergman modules.  Recall that $\mathcal M$ is said
to be a {\em quasi-free Hilbert module of rank $n$}, $1 \leq n  < \infty$, for
$\mathcal A(\Omega)$, if it is a Hilbert space completion of 
$\mathcal A(\Omega) \otimes_{\rm alg} \C^n$  (the algebraic tensor product) 
such that 
\begin{enumerate}
\item evaluation $e_w$, $e_w(f) = f(w)$, is locally uniformly bounded for 
$w \in \Omega$, 
\item pointwise multiplication by functions in $\mathcal A(\Omega)$ 
defines a bounded operator on $\mathcal M$, and 
\item  a sequence $\{f_i\}$ contained in  
$\mathcal A(\Omega)\otimes_{\rm alg}\C^n$, 
that is Cauchy in the norm of $\scM$, converges to $0$ 
in the norm of $\mathcal M$ if and only if  $\{e_w(f_i)\}$ 
converges to $0$ in $\C^n$ for all 
$w$ in $\Omega$ (cf. \cite{qfie}, \cite{quasi}).
\end{enumerate}

These assumptions ensure, among other things, via the Riesz
representation theorem, that there is a unique vector $K(\cdot, w)\in
\scM$ satisfying the reproducing property, that is,
$$h(w) = \langle h, K(\cdot, w) \rangle, ~h \in \scM,~w\in \Omega.$$

Clearly, the map $w \mapsto e_w$, which is defined on $\Omega$ and
takes values in $\scM$, is weakly holomorphic.  Hence, $e_w$ is
locally uniformly bounded in norm and $K(w,w)= \inner{e_w}{e_w}$
is locally uniformly bounded.

The Hilbert modules that we describe in this paper are the one which arise as
the quotient of a pair of Hilbert modules from the class ${\text B}_1(\Omega)$. 
The following definition makes this precise along with a mild hypothesis that 
we must impose on the quotient.   

\begin{Definition}[\cite{eqhm}, pp. 284] \label{Gen CD}
We will say that the module ${\mathcal Q}$ over the algebra
$\scA(\Omega)$ is a {\em quotient Hilbert module} in the class ${\rm B}_k(\Omega,
\scZ)$ if
\begin{enumerate} \item[(i)] there exists a resolution of the
module $\mathcal{Q}$ as in equation (\ref{quotient}) for some 
quasi-free Hilbert module $\scM$ of rank $1$ over the algebra $\scA(\Omega)$; 
\item[(ii)] for $f\in \mathcal A(\Omega)$, the restriction of the map 
$J_f$ to the hypersurface $\mathcal Z$ defines the module action on $\scQ$; 
and
\item[(iii)] the Hilbert module $J^{(k)}\scM_{|{\rm res~}\scZ}$ which is isomorphic 
to the quotient module $\scQ$ is quasi-free of rank $k$ over the algebra 
$\scA_{|{\rm res}~ \scZ}(\Omega)$. 
\end{enumerate}
\end{Definition}

Our main theorem is easily stated using the $k\times k$ array of 
differential operators:
\begin{equation} \label{diffop}
\scD_k = \begin{pmatrix}
\bar{\partial^\prime}^{\rm t} \partial^\prime  &  \bar{\partial^\prime}^{\rm t} \partial_1
&  \bar{\partial^\prime}^{\rm t} \partial_1^2&\ldots&
\bar{\partial^\prime}^{\rm t} \partial_1^{k-1} \\
\bar{\partial}_1 \partial^\prime&\bar{\partial}_1\partial_1
&\bar{\partial}_1\partial_1^2&\ldots&
\bar{\partial}_1\partial_1^{k-1} \\
\bar{\partial}_1^2\partial^\prime&\bar{\partial}_1^2\partial_1
&\bar{\partial}_1^2\partial_1^2&
\ldots& \bar{\partial}_1^2\partial_1^{k-1}\\
\vdots&\vdots&\vdots&\ddots&\vdots\\
\bar{\partial}_1^{k-1}\partial^\prime&\bar{\partial}_1^{k-1}\partial_1
&\bar{\partial}_1^{k-1}\partial_1^2&\ldots&\bar{\partial}_1^{k-1}\partial_1^{k-1}
\end{pmatrix},
\end{equation}
where ${\partial^\prime}$ denotes the differential operator
$(\partial_m, \ldots , \partial_{2})$ and $\bar{\partial^\prime}^{\rm t}$ 
is the conjugate transpose ${\partial^\prime}$. 
We point out that the $(1,1)$ position of the matrix
$\scD_k$ consists of a $(m-1) \times (m-1)$ block and that each entry
of the first row (respectively column) is a column (respectively,
row) vector of size $m-1$.

\begin{Definition} \label{}
  Let $\varrho$ and $\tilde{\varrho}$ be two positive real
  analytic functions on a domain $\Omega$. We say that the $\varrho$
  and $\tilde{\varrho}$ {\em are equivalent to order $k$ on $\scZ$} if
  $\scD_k \big ( \log \frac{\tilde{\varrho}}{\varrho} \big ) = 0$ on
  $\scZ$.
\end{Definition}

\begin{bigthm} \label{main}
Suppose that $\scQ = \scM \ominus \scM_0$ and $\tilde{\scQ} =
\tilde{\scM} \ominus \tilde{\scM_0}$ are a pair of 
quotient Hilbert modules over the algebra $\scA(\Omega)$ in the
class ${\rm B}_k(\Omega, \scZ)$. Then the quotient module $\scQ$
and $\tilde{\scQ}$ are isomorphic if and only if $\tilde{\varrho}$
and $\varrho$ are equivalent to order $k$, where $\tilde{\varrho}$
and $\varrho$ are the hermitian metrics for the line bundles
corresponding to the two modules $\scM$ and $\tilde{\scM}$
respectively.
\end{bigthm}

In case $k=2$, we may restate the Theorem \ref{main} in terms of the
tangential and the transverse curvatures along with the second
fundamental form. The validity of such a statement will follow from
Theorem \ref{main} which holds for an arbitrary $k$.  We consider the
case of $k=2$ separately for comparison with our previous result which
was limited to this case only. In the statement of the theorem below,
we use the invariants {\em tan} and {\em trans} which stand for the
tangential and transverse curvatures.  These invariants occurred in
\cite[Theorem, page 289]{eqhm}. We emphasize that the third invariant
that appears in Theorem 2 is the second fundamental form.
Therefore, this theorem is different from that of \cite{eqhm}.  We
provide an independent proof using explicit computations.

\begin{bigthm} \label{main2}
Suppose that $\scQ = \scM \ominus \scM_0$ and $\tilde{\scQ} =
\tilde{\scM} \ominus \tilde{\scM_0}$ are a pair of 
of quotient Hilbert modules over the algebra $\scA(\Omega)$ 
in the class ${\rm B}_2(\Omega, \scZ)$. Then the 
modules $\scQ$ and $\tilde{\scQ}$ are isomorphic if and
only if the restriction of the corresponding curvatures to the
hypersurface $\scZ$ coincide, that is,
\begin{description}
\item[\rm{\em tan}] $\scK_{\rm tan} = \tilde{\scK}_{\rm tan}$
\item[\rm{\em trans}] $\scK_{\rm trans} = \tilde{\scK}_{\rm trans}$
\item[\rm{\em angle}] $\scS = \tilde{\scS}$
\end{description}
are equal on $\scZ$.
\end{bigthm}

\section{Reproducing kernels and the multivariate class ${\rm B}_k$}

\subsection{} Let ${\mathcal L}(\F)$ be the Banach space of all linear
transformations on a Hilbert space $\F$ of dimension $n$ for some
$n \in \N$. Let ${\mathcal H}$ be a \hs of functions from $\Omega$
to $\F$.  For $w\in \Omega$, let $e_w:{\mathcal H}\to \F$ be
defined by $e_w(f)=f(w)$. The \hs ${\mathcal H}$ is called a
(vector-valued) functional Hilbert space if $e_w$ is bounded for
each $w\in \Omega$.  In this case, the function
$K:\Omega\times\Omega \to {\mathcal L}(\F)$ defined by
$K(z,w)=e_ze_w^*$, $z,w\in \Omega$, is called the reproducing
kernel of ${\mathcal H}$.  We recall some of the basic properties
of a reproducing kernel following \cite{Aron}.  

First, the kernel $K$ has
the reproducing property:
\begin{equation} \label{reproducing}
\inner{f}{K (\cdot,w)\eta}_{\mathcal H} =
\inner{f(w)}{\eta}_{\mathbb{F}} \mbox{~for~}\eta\in \F,~w\in \Omega,~
f\in {\mathcal H}.
\end{equation}
In particular, taking $f=K(\cdot,w)\zeta$ for $w\in\Omega$, $\zeta\in
\F$, we see that $K$ satisfies
\begin{equation} \label{pos def}
\inner{K(\cdot, w)\zeta}{K(\cdot,z)\eta} =
\inner{K(z,w)\zeta}{\eta}~~\mbox{for~}\zeta,\eta\in \F,~~z,w\in \Omega.
\end{equation}
This shows for $p \geq 1,$ $w_1,\ldots , w_p\in \Omega$ that the block
\op $\big(\!\!\big ( K(w_i,w_j)\big )\!\!\big )_{1\leq i,j \leq p}$ on
$\F\oplus\cdots \oplus \F$ ($p$ copies) is positive.
Conversely, if $K:\Omega\times \Omega \to
{\mathcal L}(\F)$ satisfies this positivity requirement for all $p$ -
tuples in $\Omega$ , one can see
that there is a unique functional \hs with reproducing
kernel $K$.  (It is the completion of the linear span of the functions
$K(\cdot, z)\eta\mbox{ for } z\in \Omega, \eta\in \F$, with inner product given
by (\ref{pos def}).)

For $1\leq i\leq m$, suppose that the \ops $M_i:{\mathcal H} \to
{\mathcal H}$ defined by $(M_i f)(z) = z_i f(z)$, $z\in \Omega, f\in
{\mathcal H}$ are bounded.  Then it is easy to verify that for each fixed
$w\in \Omega$, and $1\leq i \leq m$,
\begin{equation} \label{eigenspace}
M_i^* K(\cdot, w)\eta = \bar{w}_i K(\cdot, w)\eta \qfor \eta \in \F.
\end{equation}
Differentiating (\ref{reproducing}) we also obtain the following extension of
the reproducing property:
\begin{equation} \label{ext reprod}
\inner{(\partial_i^jf)(w)}{\eta}=\inner{f}{\bar{\partial}^j_iK(\cdot,w)
\eta}~~\mbox{for~}1\leq i\leq m,~~j\geq 0,~w\in \Omega,~\eta\in \F,~f\in
{\mathcal H}.
\end{equation}

Let $\bm{M}=(M_1, \ldots, M_m)$ denote the commuting $m$-tuple of multiplication
operators defined by the coordinate functions $z_1, \ldots , z_m$ and
$\bm{M}^*$ be the $m$-tuple $(M_1^*,\ldots ,M_m^*)$. It then
follows from (\ref{eigenspace}) that the joint eigenspace of the
$m$-tuple $\bm{M}^*$ at $w \in \Omega^*\subseteq \C^m$, where as
before, $\Omega^* = \{w\in \C^m: \bar{w} \in \Omega\}$, contains
the $n$-dimensional subspace ${\rm ran}\,K(\cdot,\bar{w})
\subseteq \scH$.

Suppose $K$ is the reproducing kernel of a \hs ${\mathcal H}$
consisting of $\F$ - valued {\em analytic functions on $\Omega$}.
Then $K$ is analytic in the first argument (and hence co-analytic
in the second argument).  We now obtain a holomorphic vector
bundle $E$ on the base space $\Omega^*$ by requiring that
$\{K(\cdot,w)v: v\in B\} \subseteq \scH$, where $B$ is an
orthonormal basis for $\F$, be a frame at $\bar{w}\in \Omega^*$.
We will also assume that $K(w,w)$ is an invertible \op for each $w\in
\Omega$.  Then $K(w,w)$ defines a hermitian metric for the bundle
$E$.  The assumption that $K(w,w)$ is invertible is automatic if 
$\mathcal H$ is a quasi-free Hilbert module of finite rank $n$.


 
Before proceeding any further, we recall the
class ${\rm B}_k(\Omega)$ for $\Omega\subset \C$ and $k\in \N$,
which was introduced in \cite{CD}.  It consists of those operators
$T$ on a Hilbert space $\scH$ for which each $w\in\Omega$ is an
eigenvalue of uniform multiplicity $k$, the eigenvectors span the
Hilbert space $\scH$ and ${\rm ran} (T-wI_{\mathcal H}) =\scH$ for
$w\in \Omega$. Later the definition was adapted to the case of an
$m$ - tuple of commuting \ops $\bm{T}=(T_1,\ldots , T_m)$ acting on a \hs $\scH$, 
first in the paper \cite{Bolyai} and then in the paper \cite{CS}
from a slightly different point of view which emphasized the role of 
the reproducing kernel.  

For $w\in \Omega \subseteq \C^m$, the $m$-tuple $\bm{T}$ is in ${\rm
B}_k(\Omega)$ if
\begin{enumerate} \label{Bk}  
\item $\rm{ran}\: D_{\mathbf{T} - w}$ is closed for all
$w\in \Omega$, where $D_{\mathbf{T}}:\scH \to \scH\otimes \C^m$
is defined by $D_{\mathbf{T}}h = (T_1h,\ldots ,T_mh)$, $h\in\scH$;
\item span~$\{\ker D_{\mathbf{T} - w}: w\in \Omega\}$
is dense in $\scH$ and
\item $\dim~\ker D_{\mathbf{T}^* - w}=k$
for all $w\in \Omega$.
\end{enumerate}
It was then shown that each of these operator $m$ - tuples $\mathbf{T}$
determines a hermitian holomorphic vector bundle $E$ of rank $k$ on
$\Omega$ and that two $m$ - tuples of operators in ${\rm B}_k(\Omega)$ are
unitarily equivalent if and only if the corresponding bundles are
locally equivalent.  In the case $k=1$, this is a question of
equivalence of hermitian holomorphic line bundles.  It is, of course,
well-known that two such line bundles are equivalent if and only if
their curvatures are equal.  However, no such simple characterization
is available if ${\rm rank}\, E =k >1$. 


We now recall that for the module $\scM$ over the algebra
$\alg{\Omega}$, the coordinate functions define an $m$-tuple of
bounded multiplication operators $\bm{M}$. We have already observed in
(\ref{eigenspace}), that the joint eigenspace of $\bm{M}^*$ at
$\bar{w},\,w\in \Omega$ includes the subspace $\{K(\cdot, w)\zeta:
w\in \Omega, \zeta \in \C^m\} \subseteq \scH$. 

Recall that if $(\scK_E)_{ij}(w) = (\scK_F)_{ij}(w)$ for $w\in\scZ$
and $2 \leq i,j \leq m$, then the restrictions of the two bundles $E$
and $F$ to the hypersurface $\scZ$ are equivalent \cite[Theorem
1.3]{rgdgm}.  In other words, $g(w) = |u(w)|^2 h(w)$ for some
holomorphic function $u$ on $\scZ$ and $w\in \scZ$.  Let $\scQ$ be the
quotient module, as in (\ref{quotient}), corresponding to the
submodule $\scM_0$ consisting of all those functions in $\scM$ which
vanish on the hypersurface $\scZ$. We showed in \cite{rgdgm} that the
restriction of just the tangential curvature $\scK_{\rm tan}$ to the
hypersurface $\scZ$ determines the quotient module up to unitary
equivalence.  If we make the stronger assumption of equality of all
coefficients of the curvature on $\scZ$, then the quotient modules can
again be shown to be equivalent in the case of $k = 2$ as is pointed
out in Remark \ref{curvrest} below. (This is a key step in the
reformulation of our earlier equivalence result for this case.)  In
this paper, we generalize this result to the case $k > 2$ by
introducing (see \eqref{diffop}) the $k \times k$ matrix $D_k$ of partial
differential operators so that the restriction of $D_k \log \varrho$
to $\mathcal Z$ determines the equivalence of the corresponding
quotient modules.  However, at this point we have no standard complex
geometric interpretation of this characterization.


\section{Jet bundles relative to a hypersurface} 
\subsection{} We are interested in submodules $\scM_0$ contained in $\scM$ 
which consist of all functions in $\scM$ that vanish to some fixed order $k$ 
on a hypersurface $\scZ$ contained in $\Omega$.  Before we can make this 
notion precise, however, we need some definitions.

\begin{Definition}[{\cite[Definition 8, p. 17]{GR}}] \label{hypersurface}

\item[(1)] A hypersurface is a complex sub-manifold of complex dimension
$m-1$, that is, a subset $\scZ \subseteq \Omega$ is a hypersurface
if for any fixed $z \in {\mathcal Z}$, there exists a
neighbourhood $U \subseteq \Omega$ of $z$ and a local defining
function $\varphi$ for $U\cap {\mathcal Z}$. 

\item[(2)] A local defining function $\varphi$ is a holomorphic map $\varphi: U
\to \C$ such that $U \cap {\mathcal Z} = \{z \in U: \varphi(z) = 0 \}$
and $\tfrac{f}{\varphi}$ is holomorphic on $U$ whenever $f_{| U\cap
  \mathcal Z} =0$. In particular, this implies that the gradient of
$\varphi$ doesn't vanish on $\mathcal Z$ and that any two defining
functions for $\scZ$ must differ by a unit.

\item[(3)] We say that the function $f$ vanishes to order $k$ on the hypersurface
$\scZ$ if $f = \varphi^n g$ for some $n\geq k$, a holomorphic function
$g$ on $U$ and a defining function $\varphi$ of $\scZ$.  The order of
vanishing on $\scZ$ of a holomorphic function $f: \Omega \to \C$ does
not depend on the choice of the local defining function.
This definition can also be framed in terms of the partial derivatives 
normal to $\scZ$.
\end{Definition}

It is clear that if there exists a global defining function
$\varphi$ for the hypersurface $\scZ$, then a valid choice of a
normal direction is the gradient of the function $\varphi$, which
is defined on all of $\Omega$.  In
general, it may not be possible to find a global defining function
for the hypersurface $\scZ$. However, if the second Cousin problem
is solvable for $\Omega$, then there exists a global defining
function (~which we will again denote by $\varphi$~) for the
hypersurface ${\mathcal Z}$.  This is pointed out in the remark
preceding Corollary 3 in \cite[p. 34]{GR}.  

Even if we don't impose the condition of ``solvability of the second
Cousin problem'' on $\Omega$, we may restrict the
holomorphic functions in the algebra $\scA(\Omega)$ and the module
$\scM$ to the open set $U$ without loss of generality. 

\subsection{} \label{z_1} We now consider in some detail the construction of 
the jet bundles needed  to characterize quotient modules of higher multiplicity.

Suppose, to begin with, we have a quasi-free Hilbert module $\scM$ of
finite rank $k$ over $\scA(\Omega)$ with kernel
function $K$ on $\Omega$ and corresponding hermitian holomorphic
vector bundle $E_{\mathcal M}$. It is easy to see that if $U$ is any
open set in $\Omega$ and $\scA_{|\rm res~U}(\Omega) = \{f_{|\rm res~U}
: f\in \scA(\Omega)\}$ is the restriction algebra, then $\scM_{|\rm
  res~U}$, the restriction of the functions in the module $\mathcal M$
to the open set $U$, is naturally a module over the algebra
$\scA_{|\rm res~U}(\Omega)$ and this module is isomorphic to $\mathcal
M$.  (Note, in general, $\scA_{|\rm res~U}(\Omega)$ is not a function
algebra since it may not be complete.)  So we can restrict all our
discussion, without loss of generality, to any open subset $U$ of
$\Omega$.  In particular, we don't need to take the domain to be as
large as possible.  Thus our treatment will be local. 


The jet bundle construction introduced in \cite{DMV} involves the kernel
function $K$ and differentiation along the normal to the hypersurface 
$\scZ \subseteq \Omega$.  We will attempt to recall the essential ideas involved as
succinctly as possible but still our description of the jet bundle will
require us to repeat substantial material from the earlier paper.

Let us construct the jet bundle $J^{(k)}E$ over some open subset $U$
of $\Omega$ which intersects $\scZ$.  Here we assume that $\mathcal M$
is rank one and make essential use of the fact that the line bundle
$E$ is given as a pullback from the Grassmanian defined by $\scM$.  In
particular, this means that a holomorphic section for $E$ over
$\Omega$ can be viewed as arising from a holomorphic function from
$\Omega$ to $\mathcal M$.  Now one takes $U$ so small that $U\cap \scZ$
equals the zero set of a holomorphic function $\varphi$ on $U$ and the
gradient of $\varphi$ doesn't vanish on $U$.

A normal direction to $\mathcal Z$ in $U \cap \scZ$ is then given by
the gradient of $\varphi$.  By choosing to reorder the coordinates and
possibly cutting down the size of $U$, we can assume that
$\frac{\partial}{\partial \!z_1} \varphi \not = 0$ on $U$.
It then follows that $\lambda_1 = \varphi(z)$,
$\lambda_2= z_2, \ldots , \lambda_m=z_m $ for $z \in \mathcal Z$
defines a local coordinate system for $U$. As pointed out
in \cite[p. 368 - 369]{DMV}, $\frac{\partial^\ell f} {\partial
z_1^\ell}(z) = 0$ for $z\in U \cap {\mathcal Z}$, $0\leq \ell \leq
k-1$ if and only if $\frac{\partial^\ell f} {\partial
\lambda_1^\ell}(\lambda) = 0$ for $\lambda\in V \cap
\phi({\mathcal Z})$, $0\leq \ell \leq k-1$, where $\phi(z) =
(\varphi(z), z_2, \ldots, z_m)$ and $V=\phi(U)$.  Then the
submodule $\scM_0$ (cf. \cite[(1.5)]{DMV}) 
consisting of those functions in $\mathcal M$ which vanish to order $k$ on the hypersurface 
$\mathcal Z$ may be described as
$$
\scM_0 = \{f\in \scM: \frac{\partial^{\ell}f}{\partial z_1^{\ell}}
(z) = 0,\:\: z\in U\cap \scZ,\:\: 0\leq \ell \leq k-1 \}.
$$
In the new coordinate
system $\phi(z) = (\varphi(z), z_2, \ldots, z_m)$, differentiation
along the normal to the hypersurface coincides with $\partial_1 =
\frac{\partial}{\partial z_1}$.  To construct the jet bundle
$J^{(k)}E$ on $U$, let us take a frame for $E$, that is, a
non-zero holomorphic section $s$ for the line bundle $E$ on $U$. 
Thus we can assume without loss of generality, 
that $s$ is defined and non-zero on all
of $U$. (Recall we can view $s$ as a holomorphic function from $U$ to $\mathcal M$.)
The jet bundle $J^{(k)}E$ over $U$ is now simply the
bundle determined by the holomorphic frame $\{s,\partial_1 s,
\ldots , \partial_1^{k-1} s\}$ on the open set $U$.
(Here, the section $s$ is viewed as a holomorphic function from $U$ to $\mathcal M$
and the differentiation of $s$ is the usual differentiation of the holomorphic function $s$.)

It is clear that the  normal direction we pick in this manner is
not unique. Thus the construction of the jet bundle $J^{(k)}E$,
even on an open subset $U$ of $\Omega$, depends on the choice of a
normal direction, and hence on the defining function.  
It follows that the normal directions obtained
from the different defining functions for $U\cap\mathcal Z$ give
rise to distinct jet bundles.  However, \cite[Proposition
2.4]{DMV} shows that these bundles coincide on $U\cap \scZ$ modulo holomorphic
hermitian equivalence. Hence, we may proceed by assuming, without loss of generality, 
that the normal direction to the hypersurface $\scZ$ is $z_1$ by
making a holomorphic change of coordinates mapping the set
$U\cap\scZ$ to $\{z_1=0\} \subseteq V \cap \Omega$ for some open
subset $V \subseteq \Omega$. Now on $U$, and using $\varphi$, we
can define a kernel function
$$JK(z,w) = (\!\!( \inner{(\partial_1^i
  s)(z)}{(\partial_1^j s)(w)}_{\mathcal M})\!\!)_{i,j=0}^{k-1}$$
on $U$. However, the kernel function $J^{(k)}K$ depends on the
choice of $U$ and $\varphi$, but the relationship between the
kernel functions obtained for different choices of sub-domains and
defining functions is particularly simple on $\scZ$.  Further, if
we restrict $JK$ to $\scZ$, or actually the intersection of this
set with $U$, then we obtain a kernel function which defines a
Hilbert space canonically isomorphic to the quotient space $\scQ$.
In Section \ref{gvsl}, we will discuss further the global versus 
local nature of the jet bundle.


\subsection{} In this subsection, we first recall the  ``change of variable formula'' 
for the jet bundle. We then define an action of the algebra $\mathcal A(\Omega)$ 
on the holomorphic sections of the jet bundle.  We use the hermitian structure of the 
jet bundle to define an inner product on the linear space of holomorphic sections of 
the jet bundle.  This is then identified as a positive definite kernel on $\Omega$.
We then discuss a notion of equivalence of the jet bundles along with a similar 
notion of equivalence for the corresponding module of holomorphic sections.   

Let us examine more closely the relationship between the jet bundle
defined by different choices of defining function on an open set $U$.
First, We recall (cf. \cite{DMV}) the construction of the jet bundle 
$J^{(k)}E$ starting with a holomorphic hermitian line bundle $E$ over
$U$. Let $s_0$ and $s_1$ be holomorphic frames for $E$ on the
coordinate patches $U_z \subseteq U$ and $U_w \subseteq U$,
respectively.  That is, $~s_0~({\rm respectively}~ s_1)~$ is a
non-vanishing holomorphic section of $E$ on $~U_z~({\rm
  respectively}~U_w)~$. Then there is a non-vanishing holomorphic
function $g$ on $U_z\cap U_w$ such that $s_0=gs_1$ there. Let
$\varepsilon_p,\:p=1,\ldots ,k$ be the standard basis vectors for
$\C^k$.  For $\ell =0,1$, we let $J(s_\ell) = \sum_{j=0}^{k-1}
\frac{\partial_1^j s_\ell}{\partial_1 w^j} \varepsilon_{j+1}$.  An
easy computation shows that $Js_0$ and $Js_1$ transform on $U_z\cap
U_w$ by the rule $J(s_0) = ({\mathcal J}g)J(s_1)$, where ${\mathcal
  J}$ is the lower triangular operator matrix
\begin{equation} \label{operator matrix}
{\mathcal J} = \begin{pmatrix}
1 & \hdots & \hdots&\hdots & \hdots& 0\\
\partial_1 & 1 & & & & \vdots \\
\vdots & & \ddots & & & \vdots \\
\vdots &\binom{l}{j}\partial_1^{\ell - j} & &1& &\vdots\\
\vdots & & & &\ddots & 0 \\
\partial_1^{k-1}&\hdots &\hdots & \hdots & \hdots& 1 \cr
\end{pmatrix}
\end{equation}
with $0 \leq \ell,j \leq k-1$.

The components of $Js$, that is, $s,\partial_1 s, \ldots ,
\partial_1^{k-1}s$, determine a frame for a rank $k$ holomorphic vector bundle
$J^{(k)}{E}$ on $U$. The transition function with
respect to this frame is represented by the matrix $({\mathcal
J}g)^{\rm tr}$, which is just the transpose of the matrix
$({\mathcal J}g)$.  We will refer to this bundle $J^{(k)}E$ over $U$ as the
$k$th order jet bundle of the bundle $E$. The
hermitian metric $\varrho(w) = \langle s(w),s(w) \rangle_{E}$ on
$E$ with respect to the frame $s$ on $E$ induces a hermitian
metric $J \varrho$ on $J^{(k)}E$ such that with respect to the
frame $Js$,
\begin{equation} \label{metric for jet bundle}
(J\varrho)(w) = \begin{pmatrix}
\varrho(w) &
\ldots &(\partial^{k-1}\varrho)(w) \cr
\vdots &
\begin{array}{ccc}
\ddots & & \\
&(\partial_1^\ell \bar{\partial}_1^m \varrho)(w)& \\
& & \ddots
\end{array}
&\vdots \cr
(\bar{\partial}_1^{k-1}\varrho)(w)& \ldots &
(\partial_1^{k-1}\bar{\partial}_1^{k-1} \varrho)(w)\cr \end{pmatrix}.
\end{equation}

Now, for any Hilbert module $\scM$ over the function algebra $\scA(\Omega)$ and
$h\in {\scM}_{| U}$, let
$$
\h = \sum_{\ell=0}^{k-1} \partial_1^{\ell} h\otimes
\varepsilon_{\ell+1}
$$
and $J(\scM_{| U}) = \{\h : h\in \scM_{| U}\}\subseteq {\scM}\otimes
\C^k$.   Consider the
map $J : \scM_{| U} \to \scM\otimes \C^k$ defined by $Jh = \h$,
for $h\in \scM_{| U}$.  Let $J^{(k)}\scM$ denote the module $J(\scM_{| U})$.
Since $J$ is injective,
we can define an inner product on $J^{(k)}\scM$
$$
\inner{J(g)}{J(h)}_{J({\scM})}=
\inner{g}{h}_{\scM}
$$
so as to make $J$ unitary. We point out that the module action on
$J^{(k)}\scM$ is no longer pointwise multiplication but the one
that ensures $J$ is a module map.

\begin{Proposition}[{\cite[page 378]{DMV}}] \label{jet kernel}
The reproducing kernel $JK: U \times U \to {\mathcal M}_k(\C)$
for the Hilbert space $J^{(k)}\scM$ is given by the formula:
$$
(JK)_{\ell,j}(z,w) =\big (\partial_1^{\ell}\bar{\partial}_1^{j}K \big )
(z,w),~~,z,w\in U,~~~ 0\leq \ell,j \leq k-1,
$$
where $\bar{\partial}_1=
\frac{\partial}{\partial\bar{w}_1}$ and $\partial_1
=\frac{\partial}{\partial z_1}$ as before.
\end{Proposition}
To complete the description of the Hilbert module $J^{(k)}\scM$, we will have to
transport the action of the algebra $\scA(\Omega)$ from $\scM$ to
$J^{(k)}\scM$ via the map $J$.  The resultant action is described in
\cite[Lemma 3.2]{DMV} which we recall now.
\begin{Lemma} \label{jet action on JD}
Let ${\mathcal M}$ be a Hilbert module of holomorphic functions on
$\Omega$ over the algebra $\alg{\Omega}$ with reproducing kernel
$K$.  Let $J^{(k)}{\mathcal M}$ be the associated module of jets with
reproducing kernel $JK$.  The adjoint of the module action  $J_f$
on $JK(\cdot,w)\x$, $\x\in \C^k$, is given by
$$
J_f^* J\kf{w} \cdot \x = J\kf{w} ({\mathcal J}f)(w)^*\cdot \x,~~f\in \scA(U),~w\in U.
$$
\end{Lemma}
The module $J^{(k)}\scM$ may be thought of as the $k$th order jet module
of the given module $\scM$ relative to the hypersurface $\scZ$.
For $w\in U$ and $1\leq \ell \leq k, $ let $s_\ell(w) = \kf{w}
\varepsilon_\ell$. The vectors $s_\ell(w)$ span the range $\scE_w$
of $\kf{w}:\C^k \to {\scM}$. The holomorphic frame $w \to
\{s_1(\bar{w}), \ldots ,s_k(\bar{w})\}$ determines a
holomorphically trivial vector bundle $\scE$ over $U^*$.  The
fiber of $\scE$ over $w$ is $\scE_w =
span\{\kf{\bar{w}}\varepsilon_\ell: 1\leq \ell \leq k\}$, $w\in
U^*$.  An arbitrary section of this bundle is of the form $s =
\sum_{\ell=1}^k a_\ell s_\ell$, where $a_\ell, \ell=1,\ldots, k,$
are holomorphic functions on $U^*$.  The norm at $w\in U^*$
is determined by
\begin{equation} \label{metric}
\|s(w)\|^2 = \inner{\sum_{\ell=1}^k a_\ell(w)
s_\ell(w)}{\sum_{\ell=1}^k a_\ell(w) s_\ell(w)}_{\scM}
=\inner{\ktf{w}{w}^{\rm tr} a(w)}{a(w)}_{{\mathbb C}^k},
\end{equation}
where $a(w) = \sum_{\ell = 1}^k a_{\ell}(w) \varepsilon_\ell$ and
$\ktf{w}{w}^{\rm tr}$ denotes the transpose of the matrix
$\ktf{w}{w}$. Since $\ktf{w}{w}$ is positive definite and
$w\mapsto \ktf{w}{w}$ is real analytic, it follows that
$\ktf{w}{w}$ determines a hermitian metric for the vector bundle
$\scE$.  It is easy to verify that if the module $\scM$ is
quasi-free and $E$ is the corresponding bundle whose hermitian
structure is determined by the kernel function $K$, then the
bundle $\scE$ along with the hermitian structure induced by the
kernel $JK$ is the one we would have obtained by applying the jet
construction to the bundle $E$.

Suppose ${\mathcal M}$ is a Hilbert module over the function
algebra $\scA(\Omega)$ which is in $B_1(\Omega)$.  Then one may
identify the Hilbert space $\scM$ with a space of holomorphic
functions on $\Omega$ possessing a complex-valued reproducing
kernel $K$.  This determines a line bundle $E_{\mathcal M}$ on
$\Omega^*$ whose fiber at $\bar{w}\in \Omega^*$ is spanned by the
vector $K(\cdot,w)$.  The jet bundle of rank $k$ is determined by
the holomorphic frame $\{\kf{w}, \bar{\partial}_1\kf{w}, \ldots,
\bar{\partial}_1^{k-1}\kf{w}\}$.  The metric for the bundle with
respect to this frame is given by the formula (compare
(\ref{metric for jet bundle}) ):
$$
\inner{\sum_{j=0}^{k-1} a_j \partial_1^j \kf{w}}{\sum_{j=0}^{k-1} a_j
\partial_1^j \kf{w}}= \sum_{j,\ell=0}^{k-1} a_j\bar{a}_\ell
\inner{\partial_1^j \kf{w}}{\partial_1^\ell \kf{w}}.
$$ 
Clearly, the action of the algebra $\scA(\Omega)$ on the module
$J^{(k)}\scM$ given in Lemma \ref{jet action on JD} defines a
holomorphic bundle map $\theta_f$ on the holomorphic frame
$\{JK(\cdot, w)\cdot \epsilon_i: 1\leq i \leq k,~w\in \Omega\}$, of
the jet bundle $J^{(k)}E_{\mathcal M}$ for each $f\in
\scA(\Omega)$. Hence the algebra $\scA(\Omega)$ acts on the
holomorphic sections of the jet bundle $J^{(k)}E_{\mathcal M}$ as well
making it into a module equivalent to the module $J^{(k)}\scM$.  This
is the jet bundle $J^{(k)}E_{\mathcal M}$ associated with $E_{\mathcal
M}$.

On the other hand, the Hilbert space $J^{(k)}{\mathcal M}$ together with its
kernel function $JK$ defined in Proposition \ref{jet kernel} defines a
rank $k$ hermitian holomorphic bundle on $\Omega^*$ (see discussion
preceding equation (\ref{metric})).  That these two constructions
yield equivalent hermitian holomorphic bundles is a consequence of the
fact that $J$ is a unitary map from ${\mathcal M}$ onto $J^{(k)}{\mathcal M}$.

\begin{Remark} \label{bundle-module}
  Therefore we see that the question of determining the equivalence
  class of the module $J^{(k)}\scM$ is the same as determining the
  equivalence class of the jet bundle $J^{(k)}E_{\mathcal M}$
  assuming that the map implementing the equivalence is also a module map
  on holomorphic sections.  Thus it is natural to make the
  following Definition.
\end{Remark}

\begin{Definition} \label{bundle-module-def}
  Two jet bundles are said to be equivalent if there is an isometric
  holomorphic bundle map which induces a module isomorphism of the
  holomorphic sections.
\end{Definition}


\section{Equivalence of jet bundles}

\subsection{} Let $E$ be a holomorphic line bundle over $\Omega^*$ equipped with
a hermitian metric $G$.  For $\scZ^* \subseteq \Omega^*$, let us
expand the real analytic function $G$ using the coordinates $(z_1,
\z) \in \Omega^*$ with $\z = (z_2, \ldots, z_m) \in \scZ^*$\,:
$$
G(z_1, \z) = \sum_{m,n =0}^\infty G_{m,n}(\z) z_1^m\bar{z}_1^n.
$$
(Note that $G$ and $G_{m,n}$ are merely real analytic functions.
Therefore, they depend on the variables we have indicated along
with their conjugates.)

Suppose we start with a resolution of the form (\ref{quotient}).  Then
we have at our disposal the domain $\Omega\subseteq \C^m$ and the
hypersurface $\scZ \subseteq \Omega$.  We recall from \cite[Theorem 3.4]{DMV}
that the quotient module $\mathcal Q$ can be identified with the 
module $J^{(k)}\mathcal M_{|{\rm res~ }\mathcal Z}$. 
The module action $J_f$ on the quotient  
$J^{(k)}\mathcal M_{|{\rm res~ }\mathcal Z}$, for 
$f\in \mathcal A(\Omega)$, is defined via the restriction of
the map 
\begin{equation} \label{nilaction}
  (J_f^* s_\ell)(w) = J\kf{w} (\mathcal Jf)(w)^* \varepsilon_\ell
\end{equation}
to $\mathcal Z$ and $\mathcal J$ is defined in equation
(\ref{operator matrix}). 

Let $\varphi$ be a local defining function for $\scZ$, that is, for
some open subset $U\subseteq \Omega$, we have $\scZ\cap U 
= \{z\in U:\varphi(z)=0\}$. If necessary, by restricting to a smaller
open subset of $U$, which we continue to denote by $U$, 
we may assume that $\varphi$ is in $\scA(U)$. 
We recall that we may assume $\Omega = U$ and $\scZ=\scZ\cap U$ without loss of
generality. Now, we see that $\varphi$ induces a nilpotent action on
each fiber of the jet bundle
$J^{(k)}E_{|{\rm res~}\scZ}$ via the restriction of the map $J_\varphi^*$
to $\scZ$.

Therefore in this picture, with the assumptions we have made along
the way, we see that the {\em quotient modules} $\scQ$ satisfies
the requirement listed in {\em (i) -- (iii)} of the following
Definition.

We begin the proof of Theorem 1 after proving a couple of results
of a general nature.  Indeed, the lemma below is a function theoretic
result and the proposition which follows is algebraic in nature.
These two results, more or less, yield immediately a proof of the theorem.

\begin{Definition}
Let $r$ be a positive real analytic function defined on $\Omega$.
Let $$ r(z_1,\z) = \sum_{\ell, m=0} ^\infty r_{\ell,m}(\z)
z_1^\ell \bar{z}_1^m
$$
be the expansion of $r$ in the variables $z_1,\, \bar{z}_1$ around
a small neighborhood of $(0,0)$. We say that {\em $r$ is
holomorphic to order $k$ along $\scZ$} if the coefficients
$r_{\ell,0},~\ell \leq k$ are holomorphic and $r_{0,m} =
\bar{r}_{m,0},~m \leq k$ are anti-holomorphic while all the
coefficients $r_{\ell,m} =0$ for $0 < \ell,m \leq k$. 
\end{Definition}

Since $\scD_2\big (\log\frac{\tilde{\varrho}}{\varrho} \big )=0$ on
all of $\Omega$ is the same as saying that 
$$
\sum_{i,j=1}^m \bar{\partial}_i\partial_j
\log\frac{\tilde{\varrho}}{\varrho} d\bar{z}_i\wedge dz_j =0 
$$
on all of $\Omega$, it follows that $\tilde{\varrho} = |\psi|^2 \varrho$ for some
holomorphic function $\psi$ on $\Omega$.  The following Lemma is a
generalization of this statement in two directions.  On the one
hand, we allow higher order differentiation and on the other hand,
we require equality only on $\scZ$.

\begin{Lemma} \label{}
Two positive real analytic functions $\varrho$ and
$\tilde{\varrho}$ on $\Omega$ are equivalent to order $k$ on
$\scZ$ if and only if $ \tilde{\varrho} = |\psi|^2  \varrho, $
where $\psi$ is some real analytic function  for which $\log
|\psi|^2$ is holomorphic to order $k$ along $\scZ$.
\end{Lemma}
\begin{proof}
  Since $\frac{\tilde{\varrho}}{\varrho}$ is a positive real analytic
  function on $\Omega$, it follows that we may write
  $\frac{\tilde{\varrho}}{\varrho} = |\psi|^2$ for some real analytic
  function $\psi: \Omega \to \C$. Let us expand the real analytic
  function $\log |\psi|^2$ in the variables $z_1$ and $\bar{z}_1$
$$
\log |\psi|^2 (z_1,\z) = \sum_{\ell, m=0} ^\infty
\psi_{\ell,m}(\z) z_1^\ell \bar{z}_1^m,
$$
where the coefficients $\psi_{\ell,m}$ are real analytic functions
of $\z\in \scZ$ for $\ell,m \geq 0$. (Strictly speaking, we should
have said $(0,\z )$ is in $\scZ$ and not $\z\in \scZ$.)

For $k=1$, to say that $\scD_1 \big ( \log
\frac{\tilde{\varrho}}{\varrho} \big ) = 0$ on $\scZ$ is the same as
saying $\partial^\prime\bar{\partial^\prime}^{\rm t}\log\frac{\tilde{\varrho}}{\varrho} =0$
on $\scZ$.  This, in turn, is equivalent to $\partial^\prime\bar{\partial^\prime}^{\rm t}
\big (\log\frac{\tilde{\varrho}}{\varrho} \big )_{|{\mathcal Z}} =0$.
As is well-known, $\partial^\prime\bar{\partial^\prime}^{\rm t} \big
(\log\frac{\tilde{\varrho}}{\varrho} \big )_{|{\mathcal Z}} =0$ if and
only if $\tilde{\varrho} = |\psi_0|^2 \varrho$ for some holomorphic
function $\psi$ on $\scZ$. This proves the Lemma for $k=1$.

The proof in the forward direction is by induction.  We have
already verified the statement for $k=1$.  Now, assume that it is
valid for $k$, that is, $\psi_{\ell,0}$ is holomorphic for $\ell
\leq k-1$, $\psi_{0,m} = \bar{\psi}_{m,0}$ for $m \leq k-1$ and
$\psi_{\ell,m} = 0$ for all $0 < \ell,m \leq k-1$.  We will show
that the same conditions are forced on the coefficients even when
we replace $k-1$ by $k$ as long as we assume $\scD_k \log
|\psi|^2=0$ on $\scZ$.  Thus we have that $\bar{\partial^\prime}^{\rm t}
\partial_1^k \log |\psi|^2_{|{\mathcal Z}} =0$ which forces
$\bar{\partial^\prime}^{\rm t}\psi_{k,0} = 0$ on $\scZ$.  Similarly,
${\partial^\prime}\psi_{0,k} = 0$ on $\scZ$ making $\psi_{0,k}$
anti-holomorphic on $\scZ$.  The condition that
$\bar{\partial}_1^\ell\partial_1^k \log|\psi|^2_{|{\mathcal Z}} =
0$ is clearly equivalent to $\psi_{\ell,k} =0$ for $\ell \leq k$.
Again, we have $\bar{\partial}_1^k\partial_1^m
\log|\psi|^2_{|{\mathcal Z}} = 0$ is clearly equivalent to
$\psi_{k,m} =0$ for $m\leq k$. Since $|\psi|^2 = \frac{1}{2}\big
(|\psi|^2 + \overline{|\psi|^2} \big )$, we see that
$\psi_{\ell,0} = \frac{1}{2} \big ( \psi_{\ell,0} +
\bar{\psi}_{0,\ell} \big )$ and $\psi_{0,\ell} = \bar{\psi}_{\ell,0}$.

The proof in the other direction is a straightforward verification --
$\scD_k \log |\psi|^2 =0$ on $\scZ$ assuming that $\psi_{\ell,0},~\ell
\leq k$ are holomorphic, $\psi_{0,m} = \bar{\psi}_{m,0},~m \leq k$ are
anti-holomorphic and the coefficients $\psi_{\ell,m} =0$ for $0 < \ell,m
\leq k$.
\end{proof}

Let $\C^{k \times k}$ be the algebra of all $k\times k$ complex
matrices and $\scT^{k\times k} \subseteq \C^{k\times k}$ be the
sub-algebra of lower triangular Toeplitz matrices, that is, those
lower triangular matrices $A$ for which $A(\ell+p,\ell) = A(p)$
for $0\leq \ell,p \leq k$, $\ell+p \leq k$.

In the proof of the following Proposition we use the fact that if
$|\psi|^2$ is holomorphic to order $k$ along $\scZ$, then the
coefficient function $\alpha_{\ell,0}$, in the expansion $|\psi|^2(\z) =
\sum_{\ell,m=0}^\infty \alpha_{\ell,m}(\z) z_1^\ell\bar{z}_1^m$,  is
a holomorphic function for $\ell \leq k$.

\begin{Proposition} \label{}
Suppose $\tilde{\varrho},\,\varrho$ are two positive real analytic functions
on $\Omega$ with $\tilde{\varrho} = |\psi|^2 \varrho$.  Then the function
$\log |\psi|^2$, which is necessarily real analytic, is holomorphic to
order $k$ along $\scZ$ if and only if there exists some holomorphic
function $\Psi: \scZ \to \scT^{k+1\times k+1}$ with $\psi_p$ at the 
$(\ell+p, \ell)$ position satisfying 
$$
(J\tilde{\varrho})(\z) = \Psi(\z) (J\varrho)(\z)
\Psi(\z)^*, ~\z\in \scZ.
$$

\end{Proposition}
\begin{proof}
Assume that $\tilde{\varrho} = |\psi|^2 \varrho$ and $|\psi|^2$
is holomorphic to order $k$ along $\scZ$.  Let us compute the derivatives
\begin{eqnarray*}
\bar{\partial}_1^i \partial_1^j \tilde{\varrho}&=&
\bar{\partial}_1^i \big (
\bar{\psi} \sum_{n_1=0}^j \binom{j}{n_2} \psi^{(n_2)} \varrho^{(j-n_2)}\big )\\
&=& \sum_{n_1=0}^i\sum_{n_1=0}^j \binom{i}{n_1} \binom{j}{n_2} \overline{\psi^{(n_1)}}
\psi^{(n_2)} \varrho^{(j-n_2, i-n_1)},
\end{eqnarray*}
where $\varrho^{(j-n_2, i-n_1)} = \partial_1^{j-n_2}\bar{\partial}_1^{i-n_1} \varrho$.
If we restrict this equation to $\scZ$, we see that
\begin{eqnarray*}
\tilde{\varrho}_{j,i} &=&  \sum_{n_1=0}^i\sum_{n_1=0}^j \bar{\psi}_{n_1} \psi_{n_2}
\varrho_{j-n_2, i-n_1},
\end{eqnarray*}
where $\tilde{\varrho}_{j,i},\, \varrho_{j-n_2, i-n_1} \mbox{ and } \psi_{n_2},\ \psi_{n_1}$
are the coefficients in the expansion of the respective real analytic functions around
$z_1^{(0)}=0$ in the variable $z_1$.  However, this says that
$
J\tilde{\varrho} = \Psi (J\varrho) \Psi^*,
$
where $\Psi$ is the lower triangular matrix with the holomorphic function $\psi_p$
at the $(\ell+p, p)$ position.

Conversely, suppose $(J\tilde{\varrho})(\z) =
\Psi_k(\z) (J\varrho)(\z) \Psi_k(\z)^*,
~\z\in \scZ,$ for some holomorphic function $\Psi:\scZ \to
\scT^{(k+1)\times (k+1)}$.  We have to show that
$\tilde{\varrho} = |\psi|^2 \varrho$ for some real analytic function
$|\psi|^2$ which is holomorphic to order $k$ along $\scZ$.
Clearly, on the hypersurface $\scZ$, we have
$$
\tilde{\varrho}_{j,i} =  \sum_{n_1=0}^i\sum_{n_2=0}^j \bar{\psi}_{n_1} \psi_{n_2}
\varrho_{j-n_2, i-n_1},
$$
where $\psi_{n_1}$ is the holomorphic function on $\scZ$ which occurs in the $n_1$
sub-diagonal of the function $\Psi$.  Now, we apply the preceding Lemma
to infer that $\tilde{\varrho}$ and $\varrho$ are equivalent to order $k$
on $\scZ$ completing the proof.
\end{proof}
\begin{Corollary} \label{jetbundleeq}
Let $(E,\varrho)$ and $(\tilde{E}, \tilde{\varrho})$ be two
hermitian holomorphic line bundles on $\Omega \subseteq \C^m$. Let
$J^{(k)}E$ and $J^{(k)}\tilde{E}$ be the jet bundles of $E$ and
$\tilde{E}$, respectively, equipped with the natural action of the
algebra $\scA(\Omega)$, that is, $f \mapsto (\scJ f) \cdot s$,
$f\in \scA(\Omega)$, for a holomorphic section $s$. The
restriction to the hypersurface $\scZ$ of the two jet bundles
$J^{(k)}E$ and $J^{(k)}\tilde{E}$ are equivalent if and only if
$\varrho$ and $\tilde{\varrho}$ are equivalent to order $k$ on
$\scZ$.
\end{Corollary}
\begin{proof}
The equivalence of the two jet bundles in the sense of Definition
\ref{bundle-module-def} amounts to the existence of a holomorphic
map $\Psi: \scZ \to \scT^{k\times k}$ which intertwines the module
action, that is,  $\Psi \scJ f = \scJ f \Psi$.  This intertwining property
is easily verified --
\begin{eqnarray*}
\Psi^* (\scJ f)^*(i,j) &=& \big ( 0,\ldots, 0, \psi_0, \ldots, \psi_{k-1-i} \big )
\big ( \partial^j f, \ldots , f, 0 ,\ldots , 0 \big )^{\mathrm tr} \\
&=& \binom{j}{i} \partial^{j-i} \psi_0+ \cdots + \psi_{j-i} f\\
&=& \psi_{j-i} f+ \cdots + \binom{j}{i} \partial^{j-i} \psi_0 \\
&=& \big ( 0,\ldots, 0, f, \ldots, \partial^{k-1-i} f \big )
\big ( \psi_j, \ldots , \psi_0, 0 ,\ldots ,0 \big )^{\mathrm tr}\\
&=&(\scJ f)^* \Psi^*(i,j)
\end{eqnarray*}
completing the proof of the Corollary.
\end{proof}

\begin{proof}[Proof of Theorem \ref{main}]
We have pointed out in Remark \ref{bundle-module} that the equivalence
of the jet bundles in the sense of Definition \ref{bundle-module-def}
is the same as that of the corresponding modules.  Therefore, the
Corollary given above completes the proof of the Theorem
\ref{main}.
\end{proof}

\section{The second fundamental form}
We let $\scM_0\subseteq \scM$ be the submodule of all functions
which vanish to order $2$ on the hypersurface $\scZ$. As before,
let $\{s, \partial_1 s\}$ be a frame for the jet bundle $J^{(2)}E$
of rank $2$ corresponding to the module $\scM$. In this case,
under some mild hypothesis on the quotient module $\scQ$, we know
\cite[pp. 289]{eqhm} that $\scK_{\rm trans}, \scK_{\rm tan}$ and
the angle $\inner{\partial s}{s}$ restricted to the hypersurface
$\scZ$ determine the unitary equivalence class of $\scQ$. Let us
explain the nature of this hypothesis.

In subsection \ref{IAS}, we show that the angle invariant, which
together with the transverse and the tangential curvatures forms a
complete set of unitary invariants for the quotient module $\scQ$ can
be replaced by the second fundamental form $\I$ for the inclusion
$E\subseteq J^{(2)}E$.  In view of the equation (\ref{F=S}), we have
stated the theorem in terms of the restriction of the curvature. One
of the disadvantages in using the angle as an invariant for the
isomorphism class of the quotient module is that for it to make sense
we must introduce normalized reproducing kernels (cf.  \cite[Remark
4.7 (b)]{CS}). To avoid this ad hoc normalization, we replace it with
the second fundamental form which is a more natural geometric
invariant.


\subsection{} \label{pati}

Let $\Omega$ be an open connected and bounded subset of $\C^m$ and
$\scZ \subseteq \Omega$ be a hypersurface, that is, a complex
sub-manifold of co-dimension $1$.  Let $\partial_1$ denote
differentiation along the normal direction to $\scZ$.
Let $E$ be a hermitian holomorphic line bundle on $\Omega$.  Let
$s$ be a holomorphic frame for $E$ and $h$ be the hermitian metric.  One sees that
$\{s, \partial_1 s\}$ is a holomorphic frame for the jet bundle $J^{(2)}E$
of rank $2$ in the normal direction to $\scZ$.
Then
\begin{equation} \label{jet2}
\big (J^{(2)} h \big ) (w)
= \begin{pmatrix} h(w) &(\partial_1 h)(w) \\
(\bar{\partial}_1 h)(w) & (\bar{\partial}_1 \partial_1 h)(w)
\end{pmatrix}
\end{equation}
defines a  metric for the jet bundle $J^{(2)}E$.
One obtains an orthonormal frame, say $\{e_1, e_2 \}$ from the
holomorphic frame by the usual Gram-Schmidt process:
\begin{eqnarray} \label{e1e2}
e_1 &=& h^{-1/2} s, \nonumber\\
e_2 &=& \frac{\partial_1 s -\inner{\partial_1 s}{e_1} e_1}
{\norm{\partial_1 s - \inner{\partial_1 s}{e_1}e_1 } } \nonumber\\
&=&\frac{\partial_1 s - \inner{\partial_1 s}{e_1} e_1}
{h^{1/2}(\partial_1\bar{\partial}_1 \log h)^{1/2}},
\end{eqnarray}
where we see that $\norm{\partial_1 s - \inner{\partial_1 s}{e_1}e_1}
= h^{1/2}(\partial_1\bar{\partial}_1 \log h)^{1/2}$ as in \cite[1.17.1]{CD}.
Let $D$ be the
canonical connection and  $\bar{\partial}$ be the operator
$\bar{\partial}f = \sum_1^m \bar{\partial_j} f d\!\bar{z}_j$.
Since $s$ is holomorphic, $\bar{\partial} s = 0$ and it follows that
\begin{equation} \label{dbare1}
\bar{\partial}e_1 = -\frac{1}{2} h^{-3/2}\bar{\partial} h \cdot s =
-\frac{1}{2}h^{-1}\bar{\partial}h  \cdot e_1 = -\frac{1}{2} \bar{\partial}(\log h) \cdot e_1.
\end{equation}
Similarly, differentiating (\ref{e1e2}), we have
\begin{eqnarray} \label{dbare2}
\bar{\partial}e_2 &=& \bar{\partial} \big ( \frac{1}
{ h^{1/2}(\partial_1\bar{\partial}_1\log h) ^{1/2}} \big )
{(\partial_1 s - \inner{\partial_1 s}{e_1} e_1)}
+ \frac{\bar{\partial} (\partial_1 s - \inner{\partial_1 s}{e_1} e_1)}
{h^{1/2}(\partial_1\bar{\partial}_1 \log h)^{1/2} } \nonumber\\
&=& - \frac{1}{2}\frac{\bar{\partial}\big ( h \partial_1 \bar{\partial}_1 \log h\big )}
{(h \deru \log h)^{3/2}} \cdot (\partial_1 s - \inner{\partial_1 s}{e_1} e_1) +
\frac{-\bar{\partial}(h^{-1} \partial_1 h )\cdot s }{(h \deru \log h)^{1/2}}  \nonumber \\
&=& - \frac{1}{2} \frac{\bar{\partial}(h (\deru \log h) )}{h \deru \log h} \cdot e_2
- \frac{\bar{\partial}(\partial_1 \log h)}{(\deru \log h)^{1/2}} e_1
\end{eqnarray}

Let us calculate the canonical hermitian holomorphic connection
$D$ in $J^{(2)}(E)$ with respect to the metric (\ref{jet2}).  We
have
\begin{eqnarray}
De_1 &=& D^{1,0} e_1 + D^{0,1}e_1 \nonumber \\
&=& \alpha_{11} e_1 +\alpha_{21} e_2 + \bar{\partial} e_1 \nonumber\\
&=& (\alpha_{11} - 1/2 \bar{\partial}\log h) e_1 + \alpha_{21} e_2
\mbox{\rm ~by (\ref{dbare1})}\nonumber\\
&=& \theta_{11} e_1 + \theta_{21} e_2,
\end{eqnarray}
where $\alpha_{11}$, $\alpha_{21}$ is a pair of $(1,0)$ forms.
Similarly, we have
\begin{eqnarray}
De_2 &=& D^{1,0} e_1 + D^{0,1}e_2 \nonumber\\
&=& \alpha_{12} e_1 +\alpha_{22} e_2 + \bar{\partial} e_2 \nonumber\\
&=&\Big ( \alpha_{12} - \frac{\bar{\partial} \partial_1 \log h}
{(\deru \log h)^{1/2}} \Big ) e_1 + \Big ( \alpha_{22} - \frac{1}{2}
\frac{\bar{\partial} \big ( h \deru \log h\big )}{h \deru \log h} \Big ) e_2
\mbox{\rm ~by (\ref{dbare2})}\nonumber\\
&=& \theta_{12} e_1 + \theta_{22} e_2,
\end{eqnarray}
where $\alpha_{12}$, $\alpha_{22}$ is another pair of $(1,0)$
forms. Since we are working with an orthonormal frame, the
compatibility with the metric (\ref{compmetric}) amounts to the
requirement
\begin{eqnarray} \label{compatible}
\{ De_i , e_j \} + \{ e_i , D e_j \} &=& \theta_{ji} + \bar{\theta}_{ij} \nonumber \\
&=& 0 \qfor 1 \leq i,j \leq 2.
\end{eqnarray}
For $1 \leq i,j \leq 2$, equating $(1,0)$ and $(0,1)$ forms separately to zero in
the equations $\theta_{ij} + \bar{\theta}_{ij} = 0$, we obtain
$\alpha_{11} = \frac{1}{2}\partial\log h$, $\alpha_{12} = 0$, $\alpha_{21} =
\frac{\partial\big ( \bar{\partial}_1 \log h \big )}
{(\partial_1 \bar{\partial}_1 \log h)^{1/2}}$ and
$\alpha_{22} = \frac{1}{2} \frac{\partial \big ( h \partial_1\bar{\partial}_1 \log h \big )}
{h \partial_1 \bar{\partial}_1 \log h}$. It therefore follows that
\begin{equation} \label{conmat}
\theta = \begin{pmatrix}
\frac{1}{2} (\partial - \bar{\partial}) \log h & - \frac{ \bar{\partial}(\partial_1 \log h)}
{(\partial_1\bar{\partial}_1 \log h)^{1/2}} \\
\frac{{\partial}(\bar{\partial}_1 \log h)}{(\partial_1\bar{\partial}_1 \log h)^{1/2}} &
\frac{1}{2} \frac{(\partial - \bar{\partial}) (h \partial_1\bar{\partial}_1 \log h )}
{h (\partial_1\bar{\partial}_1 \log h )}
\end{pmatrix}
\end{equation}
is the matrix representation of the canonical connection $D$ on
$J^{(2)}E$ with respect to the orthonormal frame $\{e_1, e_2\}$.  Thus
the second fundamental form $\I$ for the inclusion $E \subseteq
J^{(2)}E$ is
\begin{equation} \label{secfund}
\inner{De_2}{e_1} = \theta_{12} = - \frac{ \bar{\partial}(\partial_1 \log h)}
{(\partial_1\bar{\partial}_1 \log h)^{1/2}}.
\end{equation}

Let $E$ be a holomorphic hermitian vector bundle over $\Omega$. We can easily express
the second fundamental form $\I$ on $\scZ$ in terms of the coefficients of the full
curvature (\ref{curvature}) on $\scZ$:
\begin{equation} \label{F=S}
\I(z) = (\I_1(z) dz_1,\ldots, \I_m(z) dz_m) = (\scK_{\rm trans}(z))^{-1/2}
\begin{pmatrix} \scK_{\rm trans}(z) & \scS(z) \\ \end{pmatrix}
\cdot \vect{d\bar{z}_1 \wedge dz_1} {d\bar{\z}\wedge d\z}
\end{equation}
for $z = (z_1, \z) \in \Omega$.
\begin{Remark}  \label{curvrest}
It follows that if we fix the transverse curvature $\scK_{\rm trans}$ of a line bundle
$E$, then the second fundamental form $\I$ for the inclusion $E \subseteq J^{(2)}_1E$
and the coefficient $\scS$ of the curvature $\scK_E$ determine each other.
Consequently, the restriction to the hypersurface $\scZ$
of $\scK_{\rm trans}$, $\scK_{\rm tan}$ and the second fundamental form $\I$
of two holomorphic hermitian bundles are equal if and only if
the  restriction to the hypersurface $\scZ$  of all the coefficients of the
curvature $\scK$ are equal.
\end{Remark}


\subsection{} \label{IAS}

The proof of the Theorem 2 is facilitated by the following Lemma. We
let $\scK(z)$ denote the $(1,1)$ form $\sum_{i,j=1}^m
\bar{\partial}_i\partial_j (\log h)(z) d\bar{z}_i\wedge dz_j$,
$z=(z_1, \ldots , z_m) \in \Omega$, for some positive real
analytic function $h$ on the domain $\Omega$.
\begin{Lemma} \label{}
Let $h$ and $\tilde{h}$ be two positive real analytic functions on
a domain $\Omega$.  The restrictions to the hypersurface $\scZ
\subseteq \Omega$ of the corresponding $(1,1)$ forms $\scK$ and
$\tilde{\scK}$ are equal if and only if there exist holomorphic
functions $\alpha$, $\beta$ on the hypersurface $\scZ$ such that
\begin{eqnarray*}
\tilde{h}_{00} &=& |\alpha|^2 h_{00} \\
\tilde{h}_{10} &=& h_{10} + \beta |\alpha|^2 h_{00} \\
\tilde{h}_{01} &=& h_{01} + \bar{\beta} |\alpha|^2 h_{00} \\
\tilde{h}_{11} &=& |\alpha|^2 h_{11} + \beta |\alpha|^2 h_{01} +
\bar{\beta} |\alpha|^2 h_{10} + |\beta|^2 |\alpha|^2 h_{00},\\
\end{eqnarray*}
where $h(z_1, \z) = \sum_{i,j=0}^\infty h_{ij}(\z) z_1^i \bar{z}_1^j$
and $\tilde{h}(z_1, \z) =
\sum_{i,j=0}^\infty \tilde{h}_{ij}(\z) z_1^i \bar{z}_1^j$
are the power series expansions of the real analytic functions $h$ and $\tilde{h}$.
\end{Lemma}
\begin{proof}
Let us put $\gamma =\tilde{h}/h $ and
$\Gamma = \log \gamma$. Let us expand $\Gamma$ in a power series:
\begin{equation} \label{ExpandGamma}
\Gamma(z_1, \z)  = \Gamma_{00}(\z) + z_1 \Gamma_{10}(\z)+
\bar{z}_1 \Gamma_{01}(\z) + \cdots ,
\end{equation}
where $(z_1,\z)\in \Omega$. (We will suppress the dependence of the coefficients
on $\z$ whenever there is no possibility of confusion.)
Recall that $\partial^\prime = (\partial_2, \ldots , \partial_m)$.
The assumption that the restrictions to the hypersurface $\scZ \subseteq \Omega$
of $\scK$ and  $\tilde{\scK}$ are equal amounts to saying that
$\bar{\partial}\partial \Gamma = 0$. We split this condition into four separate ones.
The first of these is the requirement that
$(\bar{\partial}^\prime \partial^\prime \Gamma) _{| {\mathcal Z}} = 0$.
The second and the third
are similar: $(\partial_1 \bar{\partial}^\prime \Gamma)_{| {\mathcal Z}} = 0$
and $(\bar{\partial}_1 {\partial}^\prime \Gamma)_{| {\mathcal Z}} = 0$.  The final one is
$(\bar{\partial}_1 \partial_1 \Gamma)_{| {\mathcal Z}} = 0$.

In view of the expansion (\ref{ExpandGamma}), the first condition
is clearly the same as the requirement that $\bar{\partial}^\prime
\partial^\prime\Gamma_{00} = 0$. Therefore it follows that
$\Gamma_{00}= \alpha_1 + \bar{\alpha}_2$ for some holomorphic
functions $\alpha_1, \alpha_2$ on the hypersurface $\scZ$.  Since
$\Gamma_{00}$ is positive, we also have $\Gamma_{00}=
\bar{\alpha}_1 + {\alpha}_2$.  Hence $\Gamma_{00}=
\frac{\alpha_1+\alpha_2}{2} + \frac{\overline{\alpha_1+
\alpha_2}}{2}$.  Consequently, $\gamma_{| {\mathcal Z}} =
\exp(\Gamma_{|{\mathcal Z}}) = |\alpha|^2$, where $\alpha =
\exp{(\frac{\alpha_1 + \alpha_2}{2})}$, is a holomorphic function
defined on the hypersurface $\scZ$.

The second condition $(\partial_1 \bar{\partial}^\prime \Gamma)_{| {\mathcal Z}} = 0$
can be restated using the power series expansion (\ref{ExpandGamma}) which is
$\bar{\partial}^\prime \Gamma_{10} = 0$.  Hence $\Gamma_{10}$ is holomorphic on $\scZ$.
Similarly, $\Gamma_{01}$ is easily seen to be anti-holomorphic on $\scZ$.

Finally, the condition $(\bar{\partial}_1 \partial_1 \Gamma)_{| {\mathcal Z}}$
is clearly equivalent to the vanishing of the coefficient $\Gamma_{11}$ in the expansion
(\ref{ExpandGamma}), that is, $\Gamma_{11}=0$.

Now, we put all of the above together and modify the expansion (\ref{ExpandGamma})
\begin{equation} \label{Gamma}
\Gamma(z_1, \z) = \alpha_1 + \beta_1 z_1 + \eta_1 z_1^2+ \overline{\alpha_2 +
\beta_2 z_1 + \eta_2 z_1^2} + \cdots.
\end{equation}
It is not hard to see that we can have $\alpha_1=\alpha_2$ and
$\beta_1= \beta_2$. Indeed, $\Gamma = \frac{\Gamma +
\bar{\Gamma}}{2}$, which allows us to take the common value
$\frac{\alpha_1+\alpha_2}{2}$, and similarly
$\frac{\beta_1+\beta_2}{2}$, as the coefficient of both $z_1$ and
$\bar{z}_1$. While similar considerations apply to the coefficient
of $z_1^2$, we have to remember that in that case, and for all the
other coefficients, these are not holomorphic functions.
Therefore, we see that
\begin{eqnarray*}
\gamma &=& \exp \Gamma\\
&=& |\exp(\frac{\alpha_1+\alpha_2}{2})|^2
|\exp(\frac{\beta_1+\beta_2}{2}z_1)|^2 |\exp(\frac{\eta_1 + \eta_2}{2}z_1^2)|^2 \cdots \\
&=& |\alpha|^2 |(1+\beta z_1 + \beta^2 z_1^2 + \cdots )|^2
|(1 + \eta^2 z_1^2 + \cdots|^2 \cdots  \\
&=& |\alpha|^2 (1+ \beta z_1 +  \bar{\beta}\bar{z}_1 + |\beta|^2 \bar{z}_1z_1 + \cdots ),
\end{eqnarray*}
where $\alpha=\exp(\frac{\alpha_1+\alpha_2}{2})$ and $\beta=\frac{\beta_1+\beta_2}{2}$.
It now follows that
\begin{eqnarray*}
\tilde{h} &=& \tilde{h}_{00}+ \tilde{h}_{10} z_1 + \tilde{h}_{01} \bar{z}_1
+ \tilde{h}_{11} \bar{z}_1 z_1 + \cdots \\
&=& \gamma h \\
&=& ({h}_{00}+ {h}_{10} z_1 + {h}_{01} \bar{z}_1
+ {h}_{11} \bar{z}_1 z_1+ \cdots )(|\alpha|^2 (1+ \beta z_1 +  \bar{\beta}\bar{z}_1
+ |\beta|^2 \bar{z}_1z_1 + \cdots ))\\
&=& |\alpha|^2 (h_{00} + (h_{10} + \beta h_{00})z_1 + (h_{01} + \bar{\beta} h_{00})\bar{z}_1
+ (h_{11} + \bar{\beta} h_{10} + \beta h_{01} + h_{00} |\beta|^2) \bar{z}_1 z_1 + \cdots).
\end{eqnarray*}
Equating the coefficients in this equation, we clearly have the following relationship:
\begin{eqnarray} \label{metricrel}
\tilde{h}_{00} &=& |\alpha|^2 h_{00} \nonumber\\
\tilde{h}_{10} &=& |\alpha|^2 (\tilde{h}_{10} + \beta h_{00})\nonumber\\
\tilde{h}_{01} &=& |\alpha|^2 (\tilde{h}_{01} + \bar{\beta} h_{00})\nonumber\\
\tilde{h}_{11} &=& |\alpha|^2 (h_{11} + \bar{\beta} h_{10} + \beta h_{01} + h_{00} |\beta|^2).
\end{eqnarray}

Conversely, we see that $\scK_{22} = \frac{h_{11}h_{00} -
|h_{10}|^2}{h_{00}^2}$ on $\scZ$. If we assume the relationships
between $h$ and $\tilde{h}$ as in (\ref{metricrel}) then on the
hypersurface $\scZ$,
\begin{eqnarray*}
\tilde{\scK}_{22} &=& \frac{|\alpha|^4(h_{11} + \bar{\beta}h_{10} +\beta h_{01}
+ |\beta|^2h_{00}) h_{00} -|\alpha|^2(h_{10} +\beta h_{00})(h_{01}+\bar{\beta}h_{00})}
{|\alpha|^4 h_{00}^2} \\
&=& \frac{h_{11}h_{00} - |h_{10}|^2}{h_{00}^2}.
\end{eqnarray*}
It therefore follows that $\tilde{K}_{22} = \scK_{22}$ on $\scZ$.
Similarly, again restricted to $\scZ$, we have $\scK_{12}=
\frac{h_{00} \partial h_{01}  - h_{01} \partial h_{00}}{h_{00}^2}$.
We see that
$\scK_{12} = \frac{h_{00} \partial^\prime h_{01} - h_{01} \partial^\prime h_{00}}{h_{00}^2}$
on $\scZ$.  Hence a calculation, using (\ref{metricrel}), shows that
\begin{eqnarray*}
\tilde\scK_{12} &=& \frac{(\partial h_{01} +\bar{\beta} \partial h_{00}) h_{00} -
(h_{01} + \bar{\beta} h_{00}) \partial h_{00} }{ h_{00}^2}\\
&=& \frac{\partial^\prime h_{01}h_{00} - h_{01} \partial^\prime h_{00}}{h_{00}^2}
\end{eqnarray*}
ensuring $\tilde{\scK}_{12} =\scK_{12}$ on $\scZ$.  Finally, it is clear that
$\scK_{11}(z) = \tilde{\scK}_{11}(z)$, for $z\in \scZ$ is equivalent to
$\tilde{h}_{00} = |\alpha|^2 h_{00}$ for some holomorphic function $\alpha$
defined on $\scZ$.
\end{proof}

\begin{proof}[Proof of Theorem \ref{main2}]
We first prove the ``if'' part of the theorem.  In this case, we
have equality of all the coefficients of the two curvatures on the
hypersurface $\scZ$.  This is equivalent to the relationship given
in the equations (\ref{metricrel}).  We then find that
\begin{eqnarray*}
\begin{pmatrix} \alpha & 0\\ \alpha \beta & \alpha  \end{pmatrix}
\begin{pmatrix} h_{00} & h_{01}\\ h_{10} & h_{11}  \end{pmatrix}
\begin{pmatrix} \bar{\alpha} & \bar{\alpha} \bar{\beta}\\ 0 & 
\bar{\alpha}  \end{pmatrix}
&=&|\alpha|^2 \begin{pmatrix} h_{00} & h_{01} + \bar{\beta} h_{00}\\
h_{10} + \beta h_{00} &  h_{11} + \bar{\beta} h_{10} + \beta h_{01} 
+ h_{00} |\beta|^2
\end{pmatrix}\\
&=& \begin{pmatrix} \tilde{h}_{00} & \tilde{h}_{01}\\ 
\tilde{h}_{10} & \tilde{h}_{11}
\end{pmatrix}.
\end{eqnarray*}
It follows that the bundle map $\Theta : J^{(2)}_1E_{| {\mathcal
Z}} \to J^{(2)}_1\tilde{E}_{| {\mathcal Z}}$ defined by,
$\Theta(z) = \left ( \begin{smallmatrix} \alpha & 0\\ \alpha \beta
& \alpha
\end{smallmatrix}\right )$, for $z\in \scZ$ is holomorphic as well as isometric.
Moreover, it intertwines the nilpotent action as well. Therefore, the
quotient modules are isomorphic via this map.

For the proof of the
only if part, we first observe that any unitary implementing the
equivalence of the quotient modules must map the submodule
$\scM_0$ onto $\tilde{\scM}_0$. This implies that the tangential
curvatures must coincide. The matrix representation for the
nilpotent action corresponding to the normal coordinate has the
transverse curvature at the $(1,2)$ position.  So, if these
nilpotent actions are equivalent, then the transverse curvature
corresponding to them must coincide. Furthermore, any such
intertwining unitary between the quotient modules must be of the
form $\left (
\begin{smallmatrix} a & 0\\ b & a
\end{smallmatrix}\right )$ for holomorphic functions $a,b$ 
defined on the hypersurface $\scZ$.
We can assume, without loss of generality, that $b=a c$.  Hence we
must have
$$
\begin{pmatrix} a & 0\\ a c & a  \end{pmatrix}
\begin{pmatrix} h_{00} & h_{01}\\ h_{10} & h_{11}  \end{pmatrix}
\begin{pmatrix} \bar{a} & \bar{a} \bar{c}\\ 0 & \bar{a} \end{pmatrix} =
\begin{pmatrix} \tilde{h}_{00} & \tilde{h}_{01}\\ \tilde{h}_{10} & \tilde{h}_{11}
\end{pmatrix}.
$$
It then follows that we must have that the relationship given by
(\ref{metricrel}) holds, which completes the proof.
\end{proof}

\section{Applications and Examples}

\subsection{} Consider $\Omega_0$ contained in $\C^m$ and $\scM$ a 
quasi-free rank one 
Hilbert module for $\mathcal A(\Omega_0)$.  For $\Omega = \D \times \Omega_0$
contained in $\C^{m+1}$, we can obtain 
quasi-free rank one Hilbert modules $\mathcal R = H^2(\D)\otimes \scM$
and $\mathcal R^\prime = B^2(\D) \otimes \scM$ over $\mathcal A(\Omega)$.  
Consider the hypersurface $\mathcal Z = \{z \in \Omega : z_1 = 0\} =
0\times \Omega_0$ 
contained in $\Omega$ and the quotient Hilbert modules $\scQ = \scR/\scR_0$ 
and $\scQ^\prime = \scR^\prime/\scR^\prime_0$, where $\scR_0$ and $R^\prime_0$ 
are the submodules of functions in $\scR$ and $\scR^\prime$,
respectively, that vanish on $\scZ$. Then $\scQ \cong \scQ^\prime\cong \scM^\prime$, 
where $\scM^\prime$ is the module over $\mathcal A(\Omega)$ obtained from 
pushing forward the module $\mathcal M$ over $\mathcal A(\Omega_0)$ 
using the inclusion map $i:\Omega_0 \to \Omega$.

However, if we consider the submodules $\scR_1$ and $\scR^\prime_1$ of
functions $f$ in $\scR$ and $R^\prime$, respectively, so that both $f$
and the partial derivative of $f$ with respect to $z_1$ vanish on $\scZ$, we
obtain a rather different result.  In this case, $\scR/\scR_1 =
\scQ_1$ is not equivalent to $\scQ^\prime_1 =
\scR^\prime/\scR^\prime_1$, which can be shown by direct calculation
of the quotient modules or by using the fact that the transverse
curvatures are not equal.  

In both cases, the longitudinal curvatures agree with that of
$\scM$. In these cases, restricted to the zero set, the transverse
curvatures are constant and the angle invariant or the second
fundamental forms vanishes identically. It is not hard to produce an
example where the restriction of the transverse curvature to the zero set 
is not constant.

Let $A^2(\B^2)$ be the Bergman space on the unit ball $\B^2$.  It
consists of square integrable holomorphic functions on $\B^2$ and
possesses a reproducing kernel $B(z,w) = (1-\inner{z}{w})^{-3}$,
$z,w\in \B^2$. As it turns out, any positive real power of the Bergman
kernel $B$ is positive definite.  Therefore, there exists a Hilbert
space $A^{(\lambda)}(\B^2)$ corresponding to such a positive definite kernel
$K^{(\lambda)}(z,w):= B^{\lambda/3} (z,w) = (1-\inner{z}{w})^{-\lambda}$ 
for $\lambda >0$.  Thus we obtain a
module $A^{(\lambda)}(\B^2)$ over the polynomial algebra $\C[z]$,
$z\in \B^2$.  Now, the curvature of the corresponding holomorphic
hermitian line bundle $E^{(\lambda)}$ over the unit ball $\B^2$ is
easy to compute.  It then follows that the restriction of neither the
longitudinal nor the transverse curvature to the zero set $\{z\in \B^2:
z_1=0\}$ is constant.  However, the angle invariant is still zero
in these examples.

In the rest of this section, we construct examples of modules $\scR$ and
$\scR^\prime$ where both the longitudinal and the transverse
curvatures of these modules are the same yet the corresponding
quotient modules are not isomorphic, see Remark \ref{fund nec}.  In
these examples, it is the ``angle invariant'' which is not the same.

We also give applications of our results to a familiar
class of Hilbert modules over the bi-disc algebra.  These
applications involve homogeneity of the modules under the action
of the M\"{o}bius group. The study of homogeneity for Hilbert
modules over the algebra $\scA(\Omega)$ for a bounded symmetric
domain $\Omega\subseteq \C^m$ was initiated in \cite{MSJOT90} and
was further studied in \cite{BMJFA96}. However, in these
papers, it was assumed that $\Omega$ is irreducible.  So, the
question of considering the possibility of $\Omega=\D^2$ did not
arise.  Although, the theorem below is stated in this case, it is
clear that the proof works just as well in the case of $\D^m$. The
recent work of Ferugson and Rochberg \cite{FR} and \cite{FR1} are very close to
the discussion below -- at least, in spirit. Similarly, the work
of the second named author with Koranyi \cite{AK, KMisComptRendu} 
on homogeneous operators in the class $\mathrm B_k(\D)$ is closely 
related to what we report here.

\subsection{} For $\lambda > 0$, let $\scM^{(\lambda)}$  be the Hilbert space
which is determined by requiring that $\{e^{(\lambda)}_n(z) :=
c_n^{-1/2}z^n\,:\,n \geq 0\}$ is a complete orthonormal set in it,
where $c_n$ is the coefficient of $x^n$ in the expansion of
$(1-x)^{-\lambda}$ or $c_n$ is the set of binomial coefficients:
$\binom{-\lambda}{n} = \frac{\lambda(\lambda+1)\cdots
(\lambda+n-1)}{n!}$. It follows that $\scM^{(\lambda)}$ possesses
a reproducing kernel $K^{(\lambda)}: {\mathbb D} \times {\mathbb
D} \to {\mathbb C}$, which is given by the formula
\begin{eqnarray*}
K^{(\lambda)}(z,w) &=& \sum_{n=0}^\infty e^{(\lambda)}_n(z)
\overline{e^{(\lambda)}_n(w)}\\
&=& (1-z\bar{w})^{-\lambda},
\end{eqnarray*}
where ${\mathbb D}$ is the open unit disc.  Thus
$\scM^{(\lambda)}$ consists of holomorphic functions on the open
unit disc ${\mathbb D}$. For $\theta\in [0, 2\pi)$ and $\alpha \in
\D$, let $\varphi_{\alpha, \theta}(z) = e^{i\theta}
\frac{z-\alpha}{1-\bar{\alpha}z}$ for $z\in \D$. The group of
bi-holomorphic automorphisms M\"{o}b of the unit disc is
$\{\varphi_{\alpha, \theta}:\, \theta\in [0, 2\pi) \qand \alpha
\in \D\}$.  We recall that for $\lambda > 0$, the natural action
of the polynomial ring $\C[z]$  on each of the Hilbert spaces
$\scM^{(\lambda)}$, for $\lambda \geq 0$, makes it into a module.
However, for each $\lambda  > 1$, this action extends to the disc
algebra $\scA(\D)$. The modules $\scM^{(\lambda)}$, $\lambda \geq
0$, lie in the class $\mbox{\rm B}_1(\D)$. What is more, they are
{\em M\"{o}b -- homogeneous}, that is, the module
$\varphi_*\scM^{(\lambda)}$ defined by the action $(f,h) \mapsto
(f\circ \varphi) \cdot h$ for $f\in \scA(\D)$, $h\in
\scM^{(\lambda)}$ is isomorphic to the module $\scM^{(\lambda)}$
for all $\varphi$ in M\"{o}b.  It turns out these are the only
homogeneous modules in the class $\mbox{\rm B}_1(\D)$.  For a
complete discussion, we refer the reader to the survey paper
\cite{BMIAS01}. D. Wilkins \cite{wil} has obtained a classification of
all homogeneous Hilbert modules over the disc algebra which are in the 
class  $\mbox{\rm B}_k(\D)$ for $ k> 1$. However, he was able to
give an explicit description of these modules only for rank $2$. 
In a recent pre-print, Ferguson and Rochberg have obtained a similar
description of these modules, again only in the case of rank $2$.
A. Koranyi and the second named author have also obtained a model for 
these quotient modules \cite{AK} which works for an arbitrary $k\in
\mathbb N$.

\subsection{}
For $\lambda, \mu >0$, there is a natural action of the group
M\"{o}b$\times$M\"{o}b on the module $\scM^{(\lambda, \mu)}$,
which is just the tensor product $\scM^{(\lambda)} \otimes
\scM^{(\mu)}$.  The Hilbert space $\scM^{(\lambda, \mu)}$ is then
a space of holomorphic functions on the bi-disc via the
identification of the elementary tensor $e_m^{(\lambda)} \otimes
e_n^{(\mu)}$ with the function of two variables $z_1^m z_2^n$ on
the bi-disc ${\mathbb D} \times {\mathbb
  D}$.  It naturally possesses the reproducing kernel $K^{(\lambda,
  \mu)}(\mathbf{z}, \mathbf{w}) =
(1-z_1\bar{w_1})^{-\lambda}(1-z_2\bar{w_2})^{-\mu}$, where $\mathbf{z}
= (z_1,z_2)$ and $\mathbf{w} = (w_1,w_2)$ are both in ${\mathbb D}
\times {\mathbb D}$.
These modules are then M\"{o}b$\times$M\"{o}b -- homogeneous, with
respect to the obvious action of this group on $\scM^{(\lambda,
  \mu)}$.  We now show that these are the only M\"{o}b$\times$M\"{o}b
-- homogeneous modules which are in $\mbox{\rm B}_1(\D^2)$.

\begin{bigthm} \label{homog}
Let $\scM$ be a Hilbert module over the bi-disc algebra
$\scA(\D^2)$. Assume that $\scM$ is in ${\rm }B_1(\D^2)$ and that
it is homogeneous. Then $\scM$ is isomorphic to
$\scM^{(\lambda,\mu)}$ for some $\lambda, \mu >0$.
\end{bigthm}

\begin{proof}
Let $\gamma$ be a holomorphic section for the bundle $E$
corresponding to $\scM$.  It then follows that$\gamma\circ
\phi^{-1}$ is a holomorphic section for the module $\phi_*\scM$,
where $\phi = (\varphi_1, \varphi_2)$ is an arbitrary element of
the group M\"{o}b$\times$M\"{o}b.  These modules are then
isomorphic if and only if the curvatures of the bundle $E$
corresponding to $\scM$ and the bundle $\phi^*E$ corresponding to
$\phi_*\scM$ are equal. Let $\scK_E$ be the curvature of the line
bundle $E$, that is,
$$
\scK_E(z) =  \sum_{i,j=1}^2 \bar{\partial}_i \partial_j
 \log \|\gamma(z)\|^2 d\bar{z}_i \wedge dz_j.
$$
It will be convenient to let $\scK_E$ also denote
the coefficient matrix of the curvature of the line bundle $E$, namely
$$
\scK_E =  D \log \|\gamma\|^2, ~\mbox{where}~ D =
\big (\!\!\big ( \bar{\partial}_i \partial_j \big )\!\!\big )_{i,j=1,2}.
$$
Using the chain rule, we find that the curvature of $\phi^*E$ can
be related to the curvature of $E$ as follows.  For $\mathbf z\in
\D^2$,
\begin{eqnarray} \label{chain}
\scK_{\varphi_*E}(\mathbf z) &=& D \log \|\gamma\circ \phi^{-1}(\mathbf z)\|^2 \nonumber\\
&=& D\phi^{-1}(\mathbf z)^* \scK_E(\phi^{-1}(\mathbf z)) D\phi^{-1}(\mathbf z)\nonumber\\
&=& \begin{pmatrix}
e^{i\theta_1} \frac{1-|a_1|^2}{1+\bar{a_1}z_1} & 0\\
0 & e^{i\theta_2} \frac{1-|a_2|^2}{1+\bar{a_2}z_2} 
\end{pmatrix}^* \scK_E(\phi^{-1}(\mathbf z))
\begin{pmatrix}
e^{i\theta_1} \frac{1-|a_1|^2}{1+\bar{a_1}z_1} &0\\
0& e^{i\theta_2} \frac{1-|a_2|^2}{1+\bar{a_2}z_2} 
\end{pmatrix}. \nonumber\\
\end{eqnarray}

The equality of the curvatures for $E$ and $\phi^*E$ now
amounts to
$$
 \scK_E(\mathbf z)  =
 \begin{pmatrix}
e^{i\theta_1} \frac{1-|a_1|^2}{1+\bar{a_1}z_1} &0\\
0& e^{i\theta_2} \frac{1-|a_2|^2}{1+\bar{a_2}z_2} 
\end{pmatrix}^*
 \scK_E(\phi^{-1}(\mathbf z))
\begin{pmatrix}
e^{i\theta_1} \frac{1-|a_1|^2}{1+\bar{a_1}z_1} &0\\
0& e^{i\theta_2} \frac{1-|a_2|^2}{1+\bar{a_2}z_2}
\end{pmatrix}
$$
for all $\phi$ in M\"{o}b $\times$ M\"{o}b. By setting, $\mathbf
z=0$ in the equation relating the curvature $\scK_E$ at $\mathbf
z$ and at $\phi^{-1}(\mathbf z)$, we see that
$$
\scK_E(a_1,a_2) =
\begin{pmatrix}
e^{-i\theta_1} \frac{1}{1-|a_1|^2} &0\\
0& e^{-i\theta_2} \frac{1}{1-|a_2|^2}
\end{pmatrix}^*
 \scK_E(0,0)
\begin{pmatrix}
e^{-i\theta_1} \frac{1}{1-|a_1|^2} &0\\
0& e^{-i\theta_2} \frac{1}{1-|a_2|^2}
\end{pmatrix}.
$$
We can now put $a_1=0=a_2$ to infer that $\scK_E(0,0)$
must be diagonal, with diagonals equal to $\lambda, \mu$, say.

Finally, we can show, without loss of generality by setting
$\theta_1=0=\theta_2$, that the curvature has the form
$$
\scK_E(a_1,a_2) =  \begin{pmatrix}
\lambda (1-|a_1|^2)^{-2} & 0\\ 0 &\mu (1-|a_2|^2)^{-2}
\end{pmatrix},
$$
for $(a_1,a_2)\in \D^2$.  However, the curvature of the module
$\scM^{(\lambda, \mu)}$ has exactly this form.  So, we conclude that
the homogeneous module $\scM$ is isomorphic to $\scM^{(\lambda,
  \mu)}$.
\end{proof}

The notion of homogeneity can be adapted easily to
{\em quotient modules} over the bi-disc algebra. Let $\scM$ be a
module over the bi-disc algebra which is in the class ${\rm
B}_2(\D^2)$.  Let us define the module $\varphi_*\scM$ to be the
module which as a Hilbert space is the same as $\scM$.  However,
the algebra $\scA(\D^2)$ now acts via the map $(f,h) \mapsto
(f\circ \phi) \cdot h$, where $\phi = (\varphi,\varphi)$ with
$\varphi$ in M\"{o}b.  Let $\scM_0$ be the submodule of functions
vanishing to order $k$ on the diagonal set $\{(z,z): z\in \D\}
\subseteq \D^2$. Then the action $(f,h) \mapsto (f\circ \phi)
\cdot h$ of the algebra $\scA(\D^2)$ on $\scM$ leaves the
submodule $\scM_0$ invariant. Consequently, $\varphi_* \scM 
\subseteq \varphi_*\scM_0$.
In particular, $\varphi_* \scM^{(\lambda,\mu)}_0$ is a submodule of
$\varphi_*\scM^{(\lambda,\mu)}$. Therefore, we may form the
quotient module $\varphi_*\scQ^{(\lambda, \mu)} = \varphi_*
\scM^{(\lambda,\mu)} / \varphi_*\scM^{(\lambda,\mu)}_0$. We
clearly have, in view of Corollary \ref{qhom}, $\varphi_*\scQ$
isomorphic to $\varphi_*\scQ^{(\lambda, \mu)}$ for $\varphi$ in
M\"{o}b. This prompts the following Definition.

\begin{Definition} \label{}
Let $\exactseq{\scQ}{\scM}{\scM_0}$ be a short exact sequence of
Hilbert modules over the bi-disc algebra with the property that
the natural action of the group M\"{o}b leaves the submodule
$\scM_0$ invariant.  The quotient module $\scQ$ is said to be {\em
homogeneous} if $\varphi_*\scQ  := \varphi_*\scM / \varphi_* \scM_0$ 
is isomorphic to $\scQ$ for all $\varphi$ in the M\"{o}bius group.

\begin{Corollary} \label{qhom}
Let $\scM^{(\lambda,\mu)}_0$ be the submodule of $\scM^{(\lambda,
\mu)}$ which consist of functions vanishing to order $2$ on the
diagonal set $\triangle=\{(z,z): z\in \D\}$. Then the quotient
module $\scQ^{(\lambda,\mu)} = \scM^{(\lambda, \mu)}/\scM^{(\lambda,
  \mu)}_0$ is homogeneous.
\end{Corollary}

The proof of this Corollary is a straightforward application of the
Theorem \ref{main} which in the case of rank $2$, as we have pointed
out, says that the restriction of the curvature to the zero set is a
complete set of invariants for the quotient modules.  An explicit
description of these quotient modules follows.

\subsection{}

Let ${\scM}^{(\lambda, \mu)}_0$ be the subspace of all functions
in $\scM^{(\lambda,\mu)}$ which vanish to order $k$ on the diagonal
$\{(z,z)\,:\, z \in {\mathbb D}\} \subseteq {\mathbb D} \times {\mathbb D}$.
To describe the quotient $\scM^{(\lambda, \mu)}/\scM^{(\lambda, \mu)}_0$, it
will be useful to consider the ascending chain
\begin{equation} \label{}
 \{0\} = V_0(p) \subseteq V_1(p) \subseteq V_2(p) \cdots \subseteq V_{p+1}(p)
= \mbox{\rm Hom}(p),
\end{equation}
where Hom$(p)$ is the space of homogeneous polynomials of degree
$p$ and $V_k(p)$ is the subspace of Hom$(p)$ which is orthogonal
to the submodule $\scM^{(\lambda, \mu)}_0$.  The second named
author and B. Bagchi have developed methods to calculate
$f_p^{(k)} \in V_k(p) \ominus V_{k-1}(p)$ for $1 \leq k \leq p+1$.
These calculations are also related to the recent work of Ferguson
and Rochberg on higher order Hankel forms \cite{FR}. Also, in a recent
paper, Peng and Zhang \cite{Peng-Zhang} have shown how to carry out such 
calculations in the context of much more general domains. However, for our
purposes, we will give the details of these calculations for the
case of $k=2$ only.

First, we compute an orthonormal basis for the quotient module $\scQ =
\scM^{(\lambda, \mu)}/\scM^{(\lambda, \mu)}_0$.  We then describe the
compression of the two operators, $M_1:f\mapsto z_1f$ and $M_2: f
\mapsto z_2f$ for $f\in \scM^{(\lambda, \mu)}$, on the quotient module
$\scQ$, as a block weighted shift operator with respect to the
orthonormal basis we have computed.  These are homogeneous operators
in the class $B_2(\D)$ which were first discovered by Wilkins
\cite{wil}.

It is easily seen that
\begin{eqnarray*}
g_p^{(1)} &=& \sum_{\ell=0}^p \frac{z_1^{p-\ell} z_2^\ell}
{\|z_1^{p-\ell}\|^2 \|z_2^\ell\|^2}\\
g_p^{(2)} &=& \sum_{\ell=0}^p \frac{\ell z_1^{p-\ell} z_2^\ell}
{\|z_1^{p-\ell}\|^2 \|z_2^\ell\|^2}
\end{eqnarray*}
are in $V_1(p)$  and $V_2(p)$ respectively.  We set $f_p^{(1)} = g_p^{(1)}$.
To find $f_p^{(2)}$, all we
have to do is to find constants $a_p,\,b_p$ such that
$$\sum_{\ell=0}^p \frac{a_p\ell + b_p}{\|z_1^{p-\ell}\|^2 \|z_2^\ell\|^2} = 0.$$
This will ensure that $f_p^{(2)} = b_p g_p^{(1)} + a_p g_p^{(2)}$ vanishes on the
set $\{(z,z)\,:\,z\in \D\}$.  Hence it must be orthogonal to $V_1$.
It is clear that $a_p = - \sum_{\ell=0}^p \frac{1}{\|z_1^{p-\ell}\|^2 
\|z_2^\ell\|^2}$ and $b_p = \sum_{\ell=0}^p 
\frac{\ell}{\|z_1^{p-\ell}\|^2 \|z_2^\ell\|^2}$ meet
the requirement. Therefore,
$$
\left \{e_p^{(1)} = \frac{f_p^{(1)}}{\|f_p^{(1)}\|} , e_p^{(2)}
\stackrel{\rm def }{=} \frac{f_p^{(2)}}{\|f_p^{(2)}\|}
\right \}_{p=0}^\infty
$$
forms an orthonormal set of vectors in the quotient module 
$\scM/\scM^{(\lambda,\mu)}_2$.
To calculate the module action, we first note that
\begin{eqnarray*}
(1-|z_1|^2)^{-(\lambda+\mu)} &=&
(1-|z_1|^2)^{-\lambda} (1-|z_2|^2)^{-\mu}_{| z_1= z_2} \\
&=& \sum_{p=0}^\infty
\sum_{\ell=0}^p \frac{|z_1|^{2(p-\ell)}}{\|z_1^{p-\ell}\|^2} \frac{|z_2|^{2 \ell}}
{\|z_2^\ell\|^2}_{|z_1=z_2}\\
&=& \sum_{p=0}^\infty |z_1|^{2p} \sum_{\ell=0}^p 
\|z_1^{p-\ell}\|^{-2}\|z_2^\ell\|^{-2}.
\end{eqnarray*}
It follows that $-a_p = \|f_p^{(1)}\|^2$ is the coefficient of 
$z^p$ in the expansion of $(1-|z_1|^2)^{-(\lambda+\mu)}$ 
which is $\binom{-(\lambda+\mu)}{p}$. Similarly,
\begin{eqnarray*}
\mu (1-|z_1|^2)^{-(\lambda+\mu + 1)} &=& (1-|z_1|^2)^{-\lambda}
\frac{d}{d\, |z_2|^2}(1-|z_2|^2)^{-\mu}_{| z_1 = z_2} \\
&=& \sum_{p=0}^\infty
\sum_{\ell=0}^p \frac{|z_1|^{2(p-\ell)}}{\|z_1^{p-\ell}\|^2} 
\frac{\ell |z_2|^{2(\ell-1)}} {\|z_2^\ell\|^2}_{|z_1=z_2}\\
&=& \sum_{p=0}^\infty |z_1|^{2(p-1)} \sum_{\ell=0}^p \ell \|z_1^{p-\ell}\|^{-2}
\|z_2^\ell\|^{-2}.
\end{eqnarray*}
Therefore, we see that $b_p = \inner{g_p^{(1)}}{g_p^{(2)}}$ is the
coefficient of $z^{p-1}$ in the expansion of  $\mu
(1-|z_1|^2)^{-(\lambda+\mu+1)}$ which is $\mu
\binom{-(\lambda+\mu+1)}{p-1}$.  

Further,
\begin{eqnarray*}
\mu (1+\mu|z_1|^2) (1-|z_1|^2)^{-(\lambda+\mu + 2)} &=& (1-|z_1|^2)^{-\lambda}
\frac{d}{d\, |z_2|^2} \big (|z_2|^2 \frac{d}{d\, |z_2|^2}(1-|z_2|^2)^{-\mu}\big )
_{| z_1 = z_2} \\
&=& \sum_{p=0}^\infty
\sum_{\ell=0}^p \frac{|z_1|^{2(p-\ell)}}{\|z_1^{p-\ell}\|^2} 
\frac{\ell^2 |z_2|^{2(\ell-1)}}
{\|z_2^\ell\|^2}_{|z_1=z_2}\\
&=& \sum_{p=0}^\infty |z_1|^{2(p-1)} \sum_{\ell=0}^p \ell^2 
\|z_1^{p-\ell}\|^{-2} \|z_2^\ell\|^{-2}.
\end{eqnarray*}
Consequently, if we set $c_p = \|g_p^{(2)}\|^2$, then $c_p$ is the
coefficient of $z^{p-1}$ in the expansion of  $\mu (1+\mu|z_1|^2)
(1-|z_1|^2)^{-(\lambda+\mu + 2)} $ which is $\mu \big (
\binom{-(\lambda+\mu+2)}{p-1} + \mu
\binom{-(\lambda+\mu+2)}{p-2}\big )$. We find that
\begin{equation} \label{ac-b^2}
\|g_p^{(1)}\|^2\|g_p^{(2)}\|^2 - \inner{g_p^{(1)}}{g_p^{(2)}}^2 =
\frac{\lambda \mu}{\lambda + \mu} \binom{-(\lambda + \mu)}{p} 
\binom{-(\lambda + \mu + 2)}{p-1}.
\end{equation}
It is now easy to compute the norm of $f_p^{(2)}$:
\begin{eqnarray}
\|f_p^{(2)}\|^2 &=& \big \|\,\inner{g_p^{(1)}}{g_p^{(2)}} g_p^{(1)}-
\|g_1^{(1)}\|^2 g_p^{(2)}\, \big \|^2 \nonumber \\
&=&\|g_p^{(1)}\|^2 \big(\|g_p^{(1)}\|^2\|g_p^{(2)}\|^2 - 
\inner{g_p^{(1)}}{g_p^{(2)}}^2\big )
\nonumber \\
&=&\frac{\lambda \mu}{\lambda + \mu}\binom{-(\lambda + \mu)}{p}^2
\binom{-(\lambda + \mu + 2)}{p-1}.
\end{eqnarray}

Now, we have all the ingredients to compute the module action.
Let us first compute the matrix $M_p^{(1)} = \begin{pmatrix}
\alpha_p^{(1)} & 0\\
\beta_p^{(1)} & \eta_p^{(1)}
\end{pmatrix}$ for multiplication by $z_1$ with respect to the orthonormal basis
$\{e_p^{(1)},e_p^{(2)}\}_{p=0}^\infty$.  It is clear that
\begin{eqnarray*}
\alpha_p^{(1)} &=& \inner{z_1 e_p^{(1)}}{e_{p+1}^{(1)}}\\
&=&\frac{1}{\|g_{p+1}^{(1)}\|\,\|g_p^{(1)}\|} 
\inner{\sum_{\ell=0}^p \frac{z_1^{p+1-\ell}}
{\|z_1^{p-\ell} \|^2}
\frac{z_2^{\ell}}{\|z_2^\ell\|^2}}{\sum_{\ell=0}^{p+1} 
\frac{z_1^{p+1-\ell} z_2^\ell}
{\|z_1^{p+1-\ell}\|^2 \|z_2^\ell\|^2}}\\
&=&\frac{1}{\|g_{p+1}^{(1)}\|\, \|g_p^{(1)}\|} 
\sum_{\ell=0}^p \|z_1^{p-\ell} \|^{-2}
\|z_2^\ell\|^{-2}\\
&=& \frac{\|g_p^{(1)}\|}{\|g_{p+1}^{(1)}\|} \\
&=& \frac{\binom{-(\lambda+\mu)}{p}^{1/2}}{\binom{-(\lambda+\mu)}{p+1}^{1/2}}.
\end{eqnarray*}
Similarly,
\begin{eqnarray*}
\beta_p^{(1)} &=& \inner{z_1 e_p^{(1)}}{e_{p+1}^{(2)}} \\
&=& \frac{1}{\|g_p^{(1)}\| \|f_p^{(2)}\|} \inner{g_p^{(1)}}{f_p^{(2)}} \\
&=& \frac{1}{\|g_p^{(1)}\| \|f_p^{(2)}\|} \inner{g_{p+1}^{(1)}}{g_{p+1}^{(2)}}
\|g_p^{(1)}\|^2 - \inner{g_{p}^{(1)}}{g_{p}^{(2)}} \|g_{p+1}^{(1)}\|^2\\
&=& \big ( \frac{\mu}{\lambda}  \big )^{1/2} (\lambda + \mu +1)^{1/2} \Big (
(\lambda + \mu +p) (\lambda+\mu + p + 1) \Big )^{-1/2}\\
\end{eqnarray*}
Finally, we have
\begin{eqnarray*}
\eta_p^{(1)} &=& \inner{z_1 e_p^{(2)}}{e_{p+1}^{(2)}}\\
&=&\frac{1}{\|f_p^{(2)}\| \, \|f_{p+1}^{(2)}\|} 
\inner{z_1 f_p^{(2)}}{f_{p+1}^{(2)}}\\
&=& \frac{1}{\|f_p^{(2)}\| \, \|f_{p+1}^{(2)}\|} \|g_{p+1}^{(1)}\|^2
\big ( \|g_p^{(1)}\|^2\|g_p^{(2)}\|^2 - \inner{g_p^{(1)}}{g_p^{(2)}}^2 \big )\\
&=& \frac{\binom{-(\lambda+\mu+2)}{p-1}^{1/2}}{\binom{-(\lambda+\mu+2)}{p}^{1/2}}
\end{eqnarray*}
Since $e_p^{(2)}=0$ on the set $\{(z_1,z_2): z_1=z_2\}$, it follows that
$z_2 e_p^{(2)} = 0$ on this set as well.  Hence the projection of $z_2 e_p^{(2)}$
to the subspace $V_1(p)$ is $0$. Consequently, $M_p^{(1)}(1,2) = 0$.
Similarly, we can compute
the matrix $M_p^{(2)} = \begin{pmatrix}
\alpha_p^{(2)} & 0\\
\beta_p^{(2)} & \eta_p^{(2)}
\end{pmatrix}$ for multiplication by $z_2$ with respect to the same 
orthonormal basis
$\{e_p^{(1)},e_p^{(2)}\}_{p=0}^\infty$ as before.
Calculations similar to the ones described above show that 
$\alpha_p^{(1)} = \alpha_p^{(2)}$
and $\beta_p^{(1)} = \beta_p^{(2)}$.  However,
$\eta_p^{(2)} = - \frac{\lambda}{\mu} \eta_p^{(1)}$.

Summarizing,  the matrix
$$ M_p^{(1)} =
\begin{pmatrix}
\frac{\binom{-(\lambda+\mu)}{p}^{1/2}}{\binom{-(\lambda+\mu)}{p+1}^{1/2}} & 0\\
\big ( \frac{\mu}{\lambda}  \big )^{1/2} \frac{(\lambda + \mu
+1)^{1/2}}{\big ((\lambda + \mu +p) (\lambda+\mu + p + 1) \big
)^{1/2}} &
\frac{\binom{-(\lambda+\mu+2)}{p-1}^{1/2}}{\binom{-(\lambda+\mu+2)}{p}^{1/2}}
\end{pmatrix}$$
represents the operator $M_1$ which is multiplication by $z_1$
with respect to the orthonormal basis
$\{e_p^{(1)},e_p^{(2)}\}_{p=0}^\infty$. Similarly,
$$ M_p^{(2)} =
\begin{pmatrix}
\frac{\binom{-(\lambda+\mu)}{p}^{1/2}}{\binom{-(\lambda+\mu)}{p+1}^{1/2}} & 0\\
- \big ( \frac{\lambda}{\mu}  \big )^{1/2} \frac{(\lambda + \mu
+1)^{1/2}}{\big ((\lambda + \mu +p)(\lambda+\mu + p + 1) \big
)^{1/2}} &
\frac{\binom{-(\lambda+\mu+2)}{p-1}^{1/2}}{\binom{-(\lambda+\mu+2)}{p}^{1/2}}
\end{pmatrix}$$
represents the operator $M_2$ which is multiplication by $z_2$
with respect to the orthonormal basis
$\{e_p^{(1)},e_p^{(2)}\}_{p=0}^\infty$. Therefore, we  see that
$Q_1^{(p)} = \tfrac{1}{2}(M_1^{(p)} - M_2^{(p)})$ is a nilpotent
matrix of index $2$ while $Q_2^{(p)} =\tfrac{1}{2}(
M_1^{(p)}+M_2^{(p)})$ is a diagonal matrix in case $\mu=\lambda$.
These definitions naturally give a pair of operators $Q_1$ and
$Q_2$ on the quotient module $\scQ^{(\lambda,\mu)}$. Let $f$ be a
function in the bi-disc algebra $\scA(\D^2)$ and 
$$f(u_1, u_2)= f_{0}(u_1) + f_{1}(u_1) u_2 + f_{2}(u_1) u_2^2 + \cdots
$$
be the Taylor expansion of the function $f$ with respect to the
coordinates $u_1 = \tfrac{z_1 + z_2}{2}$ and
$u_2 = \tfrac{z_1 - z_2}{2}$. Now the module
action for $f\in \scA(\D^2)$ in the quotient module
$\scQ^{(\lambda,\mu)}$ is then given by
\begin{eqnarray*}
f \cdot h &=& f(Q_1,Q_2) \cdot h \\
& = &  f_0(Q_1) \cdot h + f_1(Q_1) Q_2 \cdot h \\
&\stackrel{\rm def}{=}& \begin{pmatrix} f_0 & 0\\ f_1 & f_0
\end{pmatrix} \cdot
\begin{pmatrix} h_1 \\ h_2 \end{pmatrix},
\end{eqnarray*}
where $h= \tbinom{h_1}{h_2} \in \scQ^{(\lambda,\mu)}$ is the
unique decomposition obtained from realizing the quotient module
as the direct sum $\scQ^{(\lambda,\mu)} = \big (\scM^{(\lambda,
\mu)} \ominus \scM^{(\lambda,\mu)}_1 \big ) \oplus \big
(\scM^{(\lambda,\mu)}_1\ominus \scM^{(\lambda,\mu)}_2 \big )$,
where $\scM^{(\lambda,\mu)}_i$, $i=1,2$, are the submodules in
$\scM^{(\lambda,\mu)}$ consisting of all functions vanishing on
$\scZ$ to order $1$ and $2$ respectively.

We now calculate the curvature $\scK^{(\lambda, \mu)}$ for the
bundle $E^{(\lambda, \mu)}$ corresponding to the metric
$K^{(\lambda, \mu)}(\mathbf{u},\mathbf{u})$, where
$\mathbf{u}=(u_1,u_2) \in \D^2$.  The
curvature $\scK^{(\lambda, \mu)}$ is easy to compute:
$$
\scK^{(\lambda, \mu)} (u_1,u_2) = (1- |u_1 + u_2|^2)^{-2}
\begin{pmatrix} \lambda & \lambda \\ \lambda & \lambda
\end{pmatrix} +(1- |u_1 - u_2|^2)^{-2}
\begin{pmatrix} \mu & - \mu \\ -\mu & \mu
\end{pmatrix}.
$$
The restriction of the curvature to the hyper-surface $\{u_2=0\}$
is
$$
\scK^{(\lambda, \mu)} (u_1,u_2)_{|u_2 = 0} = (1-|u_1|^2)^{-2}
\begin{pmatrix} \lambda + \mu & \lambda - \mu \\ \lambda - \mu &
\lambda + \mu
\end{pmatrix},
$$
where $u_1\in \D$. Thus we find that if $\lambda = \mu$, then the
curvature is of the form $2 \lambda (1-|u_1|^2)^{-2} I_2$. 

\begin{Remark} \label{fund nec} 
Let us now compare the two jet bundles, corresponding to $\lambda=\mu$
and $\lambda_1\not= \mu_1$ such that $\lambda_1 + \mu_1 =
2\lambda$. We see that the tangential and the transverse curvatures of
these line bundles  restricted to the hyper-surface
$\{u_2=0\}$ are then equal. However, the jet bundles in these two cases are not
equivalent (which is the same as saying that the quotient modules
are not equivalent). The second fundamental form, which is
"essentially" the off diagonal entry in the restriction of the curvature,
distinguishes them. In the first case it is $0$ and in the second
case it is not!
\end{Remark}

We now describe the unitary map which is basic to the construction
of the quotient module, namely, 
$$
h \mapsto {\sum_{\ell=0}^{k-1} \partial_1^{\ell} h\otimes
\varepsilon_{\ell+1}}_{|z_1=z_2}
$$
for $h \in \scM^{(\lambda, \mu)}$.  For $k=2$, it is enough to 
describe this map just for the orthonormal basis 
$\{e_p^{(1)}, e_p^{(2)}: p \geq 0\}$.  A simple calculation shows that
\begin{eqnarray}
e_p^{(1)}(z_1, z_2) &\mapsto& \begin{pmatrix}
\binom{-(\lambda+\mu)}{p}^{1/2} z_1^p \\
\mu \sqrt{\frac{p}{\lambda+\mu}} \binom{-(\lambda+\mu+1)}{p-1}^{1/2} z_1^{p-1}
\end{pmatrix} \nonumber \\
e_p^{(2)}(z_1, z_2) &\mapsto& \begin{pmatrix}
0\\
\sqrt{\frac{\lambda \mu}{\lambda+\mu}}\binom{-(\lambda+\mu+2)}{p-1}^{1/2} z_1^{p-1}
\end{pmatrix}.
\end{eqnarray}
This allows us to compute the $2\times 2$ matrix-valued kernel
function
$$
K_{\mathcal Q}({\mathbf z}, {\mathbf w}) =
\sum_{p=0}^\infty e_p^{(1)}({\mathbf z}) e_p^{(1)}({\mathbf w})^*
+ \sum_{p=0}^\infty e_p^{(2)}(\mathbf z) e_p^{(2)}(\mathbf w)^*,~{\mathbf z},
{\mathbf w}\in \D^2
$$
which restricted to $\triangle$ corresponds to the quotient Hilbert module. 
Indeed, a straight forward computation 
shows that
\begin{eqnarray*}
\lefteqn{K_{\mathcal Q}(\mathbf z, \mathbf z)_{|\mbox{res~ }\triangle} }\\
&=&\begin{pmatrix} (1-|z|^2)^{-(\lambda + \mu)} & 
\mu z (1-|z|^2)^{-(\lambda+\mu+1)}\\
\mu \bar{z} (1-|z|^2)^{-(\lambda+\mu+1)} & 
\frac{\mu^2}{\lambda+\mu} \frac{d}{d |z|^2}
\big (|z|^2(1-|z|^2)^{-(\lambda+\mu+1)}\big ) + \frac{\mu \lambda}{\lambda+\mu}
(1-|z|^2)^{-(\lambda+\mu+2)}
\end{pmatrix} \\
&=& \big (\!\! \big ( (1-|z_1|^2)^{-\lambda} {\partial^i}
\bar{\partial}^j
{(1-|z_2|^2)^{-\mu}}_{|\mbox{res~ }\triangle} \big ) \!\! \big )_{i,j =0,1}\\
&=& (JK) (\mathbf z,\mathbf z)_{|\mbox{res~ }{\mathcal Z}},\:\: \mathbf z \in \D^2,
\end{eqnarray*}
where $\triangle=\{(z,z)\in \D^2 : z \in \D\}$.
These calculations give an explicit illustration of one of the main theorems 
on quotient modules from  \cite[Theorem 3.4]{DMV}.

\subsection{}

Let  $E$ be a holomorphic hermitian line bundle defined on the bi-disc
and $J^{(k)}E$ be the jet bundle of order $k$ associated to $E$.
The M\"{o}bius group acts on the holomorphic sections of the jet
bundle $J^{(k)}E$ via the module map $s \mapsto \scJ \phi \cdot s$, where
$\phi=(\varphi , \varphi)$ for $\varphi$ in M\"{o}b.
The jet bundle along with this action of the group M\"{o}b on its
sections will be denoted by $(\scJ \varphi)^* (J^{(k)}E)$.
The bundle $E$ is said to be {\em M\"{o}b -- homogeneous of rank $k$}
if the jet bundle $J^{(k)}E$ of $E$ is equivalent to $(\scJ \varphi)^* (J^{(k)}E)$
on the set $\triangle = \{(z,z):z\in \D\} \subseteq \D^2$ for all 
$\varphi$ in the M\"{o}bius group.
\end{Definition}

It is then natural to ask which quotient modules over the bi-disc
algebra are M\"{o}b -- homogeneous. In the case of rank $k=2$, we have
shown that the modules $\scM^{(\lambda, \mu)}$ are
M\"{o}b$\times$M\"{o}b -- homogeneous.  Therefore, these are M\"{o}b --
homogeneous as well.  Are there any others?  We first consider this
question for bundles $E$ over the bi-disc.

Let $\pi: E_{\alpha,\delta}^\beta \to \D^2$ be a hermitian (trivial) holomorphic
line bundle determined by the holomorphic frame
$$
\gamma(\mathbf w)(\mathbf z) = (1-z_1\bar{w}_2)^\beta (1-z_2\bar{w}_1)^\beta
(1-z_1\bar{w}_1)^{-\alpha}(1-z_2\bar{w}_2)^{-\delta}
$$
at $\mathbf w \in \D^2$.  Let $\|\gamma(\mathbf w)\|^2 =
|(1-w_1\bar{w}_2)|^{2\beta}
(1-|w_1|^2)^{-\alpha}(1-|w_2|^2)^{-\delta}$.  We note that the
metric for the jet bundle $J^{(2)}E_{\beta}^{\alpha,\delta}$ is then
given by $\big ( \!\! \big (\partial_1^i\bar{\partial}_1^j
\|\gamma(\mathbf w)\|^2 \big ) \!\! \big )_{i,j=0,1}$. But for
this to be positive definite at $\mathbf w=(w,w)$, $w\in \D$, we
must have the conditions: $\alpha,\delta > 0$ and $\alpha\delta -
|\beta|^2 > 0$.

\begin{bigthm} \label{}
A holomorphic hermitian line bundle $E$ over the bi-disc is
M\"{o}b -- homogeneous of rank $2$ if and only if $E$ is isomorphic to
$E_{\alpha,\delta}^\beta$ for $\alpha, \delta > 0$ and some real number
$\beta$  satisfying $\alpha\delta - |\beta|^2 >0$.
\end{bigthm}
\begin{proof}
To prove the ``if'' part, we compute the curvatures of
$E_{\alpha,\delta}^\beta$ as well as that of
$\varphi^*(E_{\alpha,\delta}^\beta)$, using the chain rule
(\ref{chain}), and verify that the restrictions of these to the set
$\triangle$ are equal.

For the ``only if'' part, let $E$ be a line bundle which is M\"{o}b --
homogeneous of rank $2$.  Let $\scK({\mathbf z}) = \sum_{i,j=1}^2
\scK_{ij}({\mathbf z}) dz_i \wedge d\bar{z}_j$ be the $(1,1)$ form
valued curvature of the line bundle $E$.  Then the coefficients
${\scK_{ij}}_{| {\rm res~} \triangle}$ form a complete set of
invariants for $J^{(2)}E_{|{\rm res~} \triangle}$.

On the other hand, it is easy to see using the chain rule (\ref{chain}) that
the curvature $\scK_{\varphi^*E}$ restricted to the set $\triangle$ is given 
by the formula
$$
\scK_{\varphi^*E}(z,z) =
\frac{(1-|a|^2)^2}{|1-\bar{a}z|^4} (\scK_E \circ \phi^{-1}) (z,z),
$$
where $\phi^{-1} = (\varphi^{-1},\varphi^{-1})$ and
$\varphi(z) = \frac{z -a }{1 -\bar{a}z}$ for $a\in \D$.
Now, if $\varphi^*(J^{(2)}E)$ is unitarily equivalent to $J^{(2)}E$ on
$\triangle \subseteq \D^2$, then
$$
(\scK_E)_{ij}(z,z) = \frac{(1-|a|^2)^2}{|1-\bar{a} z|^4}
( (\scK_E)_{ij} \circ \varphi^{-1}) (z,z)
$$
for all $a\in \D$.  Putting $z=0$, we obtain
\begin{equation} \label{1}
(\scK_E)_{ij}(0,0) (1-|a|^2)^{-2} = (\scK_E)_{ij}(a,a),~ (\scK_E)_{ij}(0,0) =
\left (\begin{smallmatrix}
\alpha & \beta \\ \bar{\beta} & \delta 
\end{smallmatrix} \right ).
\end{equation}
We assume that the metric $h$ for $E$ is normalized at $0$.  The
curvature of $E$ at $0$ for a normalized metric is $\sum_{i,j=1}^2
(\bar{\partial}_i\partial_j h)(0) d\bar{z}_i \wedge dz_j$.  However,
the metric for the jet bundle $J^{(2)}E$ at $0$ is
$(\!\!(\bar{\partial}_i\partial_j h)(0))\!\!)_{i,j=1}^2$. This metric must
be positive definite which is equivalent to the condition
$\alpha\delta - |\beta|^2 > 0$.

For the rest of the proof, it will be convenient to work with the
coordinates $u_1 = (z_1+z_2)/2$ and $u_2= (z_1-z_2)/2$.
The curvature of the bundle $E$ with respect to these new
coordinates is then easily seen to be of the form
\begin{equation} \label{2}
\scK_E(u_1,u_2)_{|u_2=0} = \left(%
\begin{array}{cc}
\alpha+  \delta + \beta + \bar{\beta} &  \alpha - \delta + \beta - \bar{\beta}\\
\alpha - \delta + \beta - \bar{\beta} & \alpha + \delta -(\beta + \bar{\beta})\\
\end{array}%
\right) (1-|u_1|^2)^{-2}.
\end{equation}
Let us set $a = \alpha+  \delta + \beta + \bar{\beta}$, $b=\alpha
- \delta + \beta - \bar{\beta}$, and $c=\alpha+  \delta - (\beta +
\bar{\beta})$.  Let $\gamma(u_1,u_2) = \sum_{m,n=0}^\infty a_{mn}(u_1,
\bar{u}_1) u_2^m \bar{u}_2^n$ be a positive real analytic function on
$\D^2$. We will try to find the coefficients $a_{mn}$ so as to
ensure that the curvature of $\gamma$ restricted to the set $u_2=0$
satisfies the equation (\ref{2}).  We will let $\partial_i$ denote
differentiation with respect to $u_1$ or $u_2$ depending on whether
$i=1$ or $i=2$.  It is clear that the equation (\ref{2}) forces
\begin{eqnarray} \label{11}
\big ( \frac{\partial^{2}}{\partial_1 \bar{\partial}_1} \log
\norm{\gamma}^2\big )_{| u_2=0} &=& \frac{\partial^{2}}{\partial_1
\bar{\partial}_1} \log a_{00} \\ \nonumber 
&=& a (1-|u_1|^2)^{-2}.
\end{eqnarray}
It then follows that $a_{00} = (1-|u_1|^2)^{-a}$.  Similar calculations show that
\begin{eqnarray} \label{12}
\big ( \frac{\partial^{2}}{\bar{\partial_1}{\partial}_2} \log
\norm{\gamma}^2\big )_{| u_2=0} &=& a_{00}^{-2} (a_{00}\,
\bar{\partial}_1\, a_{10} - a_{10}\, \bar{\partial}_1 \,a_{00}) \\ \nonumber
&=& b (1-|u_1|^2)^{-2}.
\end{eqnarray}
Choosing $a_{10} = (b/a) \partial_1\, a_{00} = b\, \bar{u}_1 (1-|u_1|^2)^{-a-1}
= b\, a_{00} \,\bar{u}_1(1-|u_1|^2)^{-1} $, we verify the equation (\ref{12}).
Finally, we have
\begin{eqnarray} \label{22}
\big ( \frac{\partial^{2}}{\bar{\partial_2}{\partial}_2} \log
\norm{\gamma}^2\big ) _{| u_2=0} &=& a_{00}^{-2} ( a_{11} a_{00} -
|a_{10}|^2) \\ \nonumber 
&=& c (1-|u_1|^2)^{-2}.
\end{eqnarray}
We can now solve for
\begin{eqnarray} \label{22s}
a_{11} &=& a_{00}^{-1} \big ( c(1-|u_1|^2)^{-2} a_{00}^2 + b^2 |u_1|^2 a_{00}^2 
(1-|u_1|^2)^{-2} \big ) \\ \nonumber
&=& a_{00} (1-|u_1|^2)^{-2} (c + b^2|u_1|^2).
\end{eqnarray}
 
Recall that the restriction of the curvature
determines the coefficients $a_{00}, a_{10}$ and $a_{11}$ in the
metric $\gamma$ modulo unitary equivalence of the quotient
modules.  Therefore, the positive definite matrix-valued function
$\Gamma = \begin{pmatrix} a_{00} & a_{01} \\ a_{10} & a_{11}
\end{pmatrix}$ describes all possible {\em homogeneous} bundles of rank $2$
on the bi-disc.  We see that the jet bundle  
$J^{(2)}E_{\alpha,\delta}^\beta$ on $\triangle$ 
may be obtained  from the line bundle
$E_{\alpha,\delta}^\beta$ and that the curvature of  this line bundle,
computed with respect to the variables $u_1,u_2$ at $(u_1,0)$, 
$u_1\in \D$, is exactly what is prescribed in (\ref{2}). This
completes the proof. \end{proof}

Whether the holomorphic hermitian line bundles
$E_{\alpha,\delta}^\beta$, $\alpha\delta -|\beta|^2>0$, correspond to
a Hilbert module $\scM$ over the algebra $\scA(\Omega)$ depends on the
question of positive definiteness of the function $\gamma(\mathbf
w)(\mathbf z)$ for $\mathbf z,\,\w\in \D^2$.

\section{Some closing Remarks} \label{gvsl}
As is true in many cases, the current paper probably raises as many
questions as it answers.  While our hope is to investigate many of the
directions suggested in the future, we want to point them out here.
Also, other thoughts seem to be of a more intuitive, preliminary nature
but promise tantalizing connections with other topics.  We will attempt
to record these possibilities as well.
\subsection{} We begin with a succinct conceptual recollection of the original
connection of operator theory with complex geometry couched in the
context studied in this paper.

As mentioned in Section 1, the kernel function $K_{\mathcal M}$
defined for a finite-rank $k$ quasi-free Hilbert module $\mathcal M$
over a domain $\Omega$ can be used to define a hermitian holomorphic
rank $k$ vector bundle $E_{\mathcal M}$ which is a pullback of a
holomorphic map from $\Omega$ to the Grassmanian of $k$-dimensional
subspaces of $\mathcal M$.  Moreover, this bundle determines the
module up to unitary equivalence.  Since for $U$ an open subset of
$\Omega$, one can show that the span of the fibers of $E_{\mathcal M}$
over $U$ equals $\mathcal M$, the restriction $(E_{\mathcal M})_{|U}$
of $E_{\mathcal M}$ to $U$ also determines $\mathcal M$.  Hence there
is no compelling reason to consider the bundle over the largest open
set possible.  However, the fibers of $E_{\mathcal M}$ over any point
of $\Omega$ can still be seen in terms of $\mathcal M$.

In particular, the fiber of $E_{\mathcal M}$ at $w \in \Omega$ can be
identified naturally with the quotient $\scM/[\scA(\Omega)_w \scM]$,
where $[\scA(\Omega)_w \scM]$ denotes the closure of the linear span
of the products of $\scA(\Omega)_w$ with the functions in $\scM$, and
$\scA(\Omega)_w$ is the maximal ideal of functions in $\scA(\Omega)$
that vanish at $w$.  It is shown in \cite{DMV} that the disjoint union
of these fibers can be identified with $E_{\mathcal M}$.  Moreover,
for $f$ a function in $\scA(\Omega)$, the module action defines a
holomorphic bundle map on $E_{\mathcal M}$ which is multiplication by
the scalar $f(w)$.  We complete this brief summary by stating the
three basic parts of the theory in the form of a theorem.

\begin{Theorem} \label{7.1}
Let $\scM$ be a Hilbert module in the class ${\rm B}_k(\Omega)$ with 
associated bundle $E_{\mathcal M}$.  Then
\begin{enumerate}
\item[(a)]  a complete set of ``geometric invariants'' for a hermitian  
holomorphic vector bundle $E$, which determines the bundle up to 
equivalence, consists of its curvature and sufficiently many 
partial derivatives of the curvature;
\item[(b)] a complete set of ``operator invariants'' which determines  
the Hilbert module $\scM$ up to unitary equivalence, consists of the 
$m$-tuples of commuting nilpotent matrices obtained by restricting 
the coordinate multiplication operators to the common generalized 
eigenspaces to high enough order; and
\item[(c)]  the ``geometric invariants''  of (a) determines the ``operator 
invariants'' of (b) and vice versa.
\end{enumerate}
\end{Theorem}

We refer the reader to the earlier papers \cite{CD, Bolyai, CD85, chenrgd} 
for complete details. 

\subsection{} 
Now we want to consider the same set of questions for the quotient
Hilbert modules considered in this paper.
 
We begin with a few comments on the notion of an analytic
hypersurface.  In general, a subset $\scZ$ of $\Omega$ defined as the
zero set of a holomorphic function possesses singularities of various
kinds.  Even so, the set of smooth manifold points forms a dense open
subset $\scZ^\prime$ of $\scZ$.  Although one can restrict attention
to $\scZ^\prime$ contained in a smaller open subset of $\Omega$, as we
have done, a function in $\scM$ that vanishes on $\scZ^\prime$ will
actually vanish on all of $\scZ$.  Moreover, the quotient Hilbert
module will yield a kind of spectral sheaf defined over \textit{all} of
$\scZ$ with the fibers over singular points also having an operator
theoretic meaning.  But this phenomenon is a topic for a later
investigation.  Thus we will assume, as we have done in the paper,
that $\scZ$ is a smooth manifold.
 
In Section 3 we showed how to construct the jet bundle $J^{(k)}
E_{\mathcal M}$ over an open subset $U$ of $\Omega$ on which there is
a ``good defining function'' $\varphi$ and determined the change in
this construction corresponding to a change in defining function.  An
obvious question which presents itself at this point is whether or not
the jet bundle can be defined over a neighborhood of $\scZ$ or even on
all of $\scZ$?  As in the previous section, one can use the fact that
the $J^{(k)} E_{\mathcal M}$ constructed on an open set $U$ is defined
as a pullback bundle from the Grassmanian of $k$-dimensional subspaces
of the quotient Hilbert space $\scQ$ to identify it and its fibers
concretely, at least over points of $\scZ$.  Analogous to the earlier
case, such a fiber can be identified   with $\scQ/[\scA^{(k)}(\Omega)_w
\scQ]$ for $w$ in $\scZ$.  Thus, one can show that $J^{(k)} E_{\mathcal M}$ 
is a well-defined hermitian holomorphic vector bundle over $\scZ$.  
Actually, there is a simpler expression for
these fibers.  For $w \in \scZ$ and $v$ a vector normal to $\scZ$ at
$w$, let $\scA(\Omega)_{w,v}$ denote the functions in $\scA(\Omega)$
for which both the function and the partial derivative in the
$v$-direction vanish at $w$ to order $k$.  Then one can show that 
$\scQ/[\mathcal A^{(k)}(\Omega)_w \scQ]$ is naturally isomorphic to 
$\scM/[\scA^{(k)}(\Omega)_{w,v} \scM]$.  (Here the exponent again 
refers to the linear span of $k$-fold products.) 
In this context, even more is true.  
 
The identification of $\scQ/[\mathcal A^{(k)}(\Omega)_{w,v}\scQ]$
  with the fiber over $w$ preserves $\scQ/[\scA^{(i)}(\Omega)_{w,v}
  \scQ]$ for $1 \leq i \leq k$, and hence the flag structure of
  $J^{(k)} E_{\mathcal M}$ is also well-defined over $\scZ$.  To make
  this more precise, one needs to recall the special frame for the jet
  bundle constructed over an open set $U$ in Section 4.  Now the
  metric on $J^{(k)} E_{\mathcal M}$ defined in Section 3 is the same
  as the one inherited from the Grassmanian or the quotient norm on
  $\scM/[\scA^{(k)}(\Omega)_{w,v} \scM]$.  But there is even more
  structure present.
 
For $\psi$ a function in $\scA(\Omega)$, a bounded operator is defined
on $\scQ$ and hence also on each fiber $\scQ/[\scA^{(k)}(\Omega)_{w,v} \scQ]$.
Relative to the special basis chosen in Section 4, the operator at
each point $w$ is a Toeplitz-like matrix.  In particular, the matrix for a
defining function for $\scZ$ at $w$ in $\scZ$ is a nilpotent matrix of
order $k$.  It is the unitary equivalence class of this nilpotent at
the points $w$ in $\scZ$ that corresponds to the operator invariants for
this case. We summarize these results in the following theorem: 
\begin{Theorem} \label{7.2}
Let $\scM$ be a rank one quasi-free Hilbert module over $\scA(\Omega)$
and $\scZ$ be an analytic hypersurface contained in $\Omega$.  Then
the jet bundle $J^{(k)} E_{\mathcal M}$ over $\scZ$ can be identified
with the union of the fibers $\scM/[\scA^{(k)}(\Omega)_{w,v} \scM]$.
Moreover, the module action induces by restriction to each fiber an
algebra isomorphic to the lower triangular Toeplitz matrices.
Finally, the quotient module determines these fiber operators up to
unitary equivalence.
\end{Theorem}

Unfortunately, at this point we don't understand what constitutes a
complete set of ``operator-theoretic invariants'',
although, in analogy with the results described in Section 7.1, we
might expect it to be the commuting $m$ - tuple of nilpotents obtained
from the restriction of the coordinate multipliers to higher order
generalized eigenspaces.  These invariants can be viewed as analogues
to ``geometric invariants'' but except for the case $k = 2$, a better
description should be possible.  We will say more about this matter
below.  We were, however, able to obtain a complete set of invariants
in terms of the operator $D_k$ which is the result we presented in
Section 4.

\subsection{} 
In this part we begin by reviewing what it means for bundles to be
equivalent in terms of frames, in both the contexts of sections 7.1 and 7.2.
With that information in hand, we will see that characterizing
equivalence can be divided into two parts, equivalence at a point and
equivalence in a neighborhood of the point.  After that, we will
attempt to use this framework to interpret the invariants we have
obtained earlier in the paper.

In section 7.1 the bundle $E_{\mathcal M}$ in question has rank $k$
and a hermitian
holomorphic structure and is defined as a pullback from the
Grassmanian.  Moreover, at least locally on an open set $U$ of $\Omega$,
one can find a holomorphic frame $\{s_1 (w),\ldots , s_k (w)\}$ which we can
take to be holomorphic $\mathcal M$ - valued functions on $U$, where 
$\mathcal M$ is a Hilbert module over $\scA(\Omega)$.

Now suppose $\tilde{\mathcal M}$ is another Hilbert module over
$\scA(\Omega)$ which defines a rank $k$ hermitian holomorphic bundle
$E_{\tilde{\mathcal M}}$ with a holomorphic frame
$\{\tilde{s}_1(w),\ldots , \tilde{s}_k(w)\}$ also over $U$.  What does
it mean to say that $E_{\mathcal M}$ and $E_{\tilde{\mathcal M}}$ are
equivalent over $U$?

Essentially, there must exist a $k \times k$ matrix of holomorphic functions
$(\!\!(\psi_{i,j})\!\!)$ on $U$ such that
\begin{enumerate}
\item[(1)] $\tilde{s}_p(w) = \sum_{j=1}^k \psi_{p,j}(w) s_j(w)$ for $p =
1,\ldots, k$; and
\item[(2)] the matrix $(\!\!( \psi_{i,j}(w))\!\!)$ defines a unitary
map between the corresponding fibers of $E_{\mathcal M}$ and
$E_{\tilde{\mathcal M}}$ for $w \in U$.
\end{enumerate}
Now an obvious necessary condition for the existence of such a matrix
of functions is that such a matrix must exist at each point $w$.  This
is the pointwise condition mentioned above.  However, here that
condition is vacuously satisfied.

By hypothesis, the set of values at $w$ of a frame over $U$ for
$E_{\mathcal M}$ forms a basis for the $k$-dimensional fiber as does
the set of values at $w$ of a frame over $U$ for $E_{\tilde{\mathcal
M}}$.  Now both fibers have an inner product and we can find a matrix
taking one basis to the other and acting as a unitary. (Note this is
not the same thing as saying that the matrix is a unitary matrix since
the inner products on the domain and range are different.)  However,
note that such a matrix is far from being unique since we can both
pre- and post-multiply it by a unitary matrix.  In case $k = 1$, or
the bundles are line bundles, the non-uniqueness is a scalar of
modulus one.

In the general case, one can choose a matrix of functions which
accomplishes both (1) and (2) but the question is whether or not those
functions can be chosen to be holomorphic.  That is the question
answered in \cite{CD}, \cite{Bolyai}, and \cite{CD85} with the answer
involving the curvature and partial derivatives of the curvature.  We
will not proceed any further with a descriptive analysis in this case.

\subsection{} 
Now we want to treat the bundle discussed in 7.2 which arises from the
quotient Hilbert module in the same fashion as we did for $E_{\mathcal
M}$.

In particular, we have the jet bundles for the two Hilbert modules
$\mathcal M$ and $\tilde{\mathcal M}$.  Each has rank $k$ and there is
a canonical frame $s(w)$ over an open set $U$ for each once one fixes
sections $s(w)$ and $\tilde{s}(w)$.  The other elements of the frame
are obtained by differentiating the given section in the direction
normal to the hypersurface $\scZ$ using the same good defining
function for each.  Again, we ask when these two bundles are
equivalent but now we want more, not just equivalence of the two
bundles but a bundle map effecting that equivalence which is also a
module map.  Before discussing just what that entails, let us point
out that although we didn't mention it in 7.3, the bundle maps
discussed there were module maps because the action induced by a
multiplier $\psi$ in $\scA(\Omega)$ on the fiber over $w$ is just
multiplication by the scalar $\psi(w)$.

Now a bundle map effecting an equivalence between $J^{(k)}E_{\mathcal
M}$ and $J^{(k)} E_{\tilde{\mathcal M}}$ must again be a matrix of
holomorphic functions satisfying (1) and (2) but now there is also a
condition:
\begin{enumerate}
\item[(3)] the matrix for the value of $(\!\!(\psi_{i,j})\!\!)$ at $w$ 
is a Toeplitz-like matrix, that is, it is lower triangular and the entries 
on a diagonal are predetermined multiples of each other.  Moreover, 
the matrix corresponding to a defining function
for $\scZ$ at $w$ is a nilpotent of order $k$ which is a single Jordan block.
\end{enumerate}
This latter condition places strong restrictions on the matrix
function $\Psi$, particularly in view of (2) which
means it must define a unitary map.  And, whereas in the case of 7.3
there is no pointwise obstruction, now there is.  This issue can be
approached as follows.

Consider a separable, infinite dimensional Hilbert space $\mathcal H$
and the collection $N_k(\scH)$ of ordered, linearly independent
$k$-tuples $X = \{x_1, \ldots, x_k\}$ in $\scH$.  For a given $X$ in
$N_k(\scH)$, there is a unique element ${\rm Gr}(X)$ of the
Grassmanian, ${\rm Gr}_k(H)$, of $k$ - dimensional subspaces of $\scH$
which it determines.  There is also an element ${\rm St}(X)$ in the
complex Stiefel manifold of linearly independent subsets with $k$
elements.  Finally, let us consider the order $k$ nilpotent operator
${\rm Nil}(X)$ defined on the span of the vectors in $X$ by the simple
shift, that is, the operator which takes $x_i$ to $x_{i+1}$ for $0
\leq i < k$.

We can define several notions of equivalence on $N_k(\scH)$ as
follows.  First, we can identify $X$ and $X^\prime$ if the subspaces
they span are equal or equivalently, if ${\rm Gr}(X) = {\rm
  Gr}(X^\prime)$.  Second, we can identify them if the two Stieffel
elements, ${\rm St}(X)$ and ${\rm St}(X^\prime)$, are unitarily
equivalent. Finally, we can identify them if the nilpotent operators
${\rm Nil}(X)$ and ${\rm Nil}(X^\prime)$ are unitarily equivalent.  One can
easily see that equivalence of the nilpotents implies equivalence of
the Stieffel elements which in turn implies equivalence of the Grassmanians, and
none of the equivalences are the same.  Moreover, one can easily
determine the Lie group of operators that respect each of the
equivalences, in case $X = X^\prime$.

Now let us return to the question of a pointwise obstruction to the
existence of a $k \times k$ matrix of holomorphic functions satisfying (1),
(2) and (3).  

\begin{Theorem} \label{7.3}
        Let $s(w)$ and $s'(w)$ be the canonical  frames over $U$ for
two jet bundles determined by the same defining function and consider
the elements $S(w)$ and $S'(w)$ of $N_k(\scM)$ and $N_k(\scM^\prime)$,
respectively, that they determine.  Then a necessary condition for the
jet bundles to be equivalent is that $N_k(\scM)$ and
$N_k(\scM^\prime)$ are equivalent.
\end{Theorem}

The proof is straightforward since (1), (2), and (3) imply equivalence 
of the nilpotents. 

If $s(w)$ and $s^\prime(w)$ are the canonical frames over $U$ for the
two jet bundles determined by the same defining function, then for
each $w \in U$ they yield elements in $N_k(\mathcal M)$ and
$N_k(\mathcal M^\prime)$, respectively, by evaluating the ordered
frames at $w$.  Conditions (1), (2) and (3) imply that the
corresponding nilpotent operators are unitarily equivalent.  Hence for
each $w$, a necessary condition for the jet bundles to be equivalent
is that the elements in $N_k$ are equivalent, and this does not always
happen.  The relationship of this condition to the unitary invariants
obtained in this paper will be considered in subsequent work.

\subsection{} We conclude with a number of comments suggesting additional
connections or further lines of investigation of the results of this
paper.

The ``nilpotent invariants'' identified in the previous subsection
refine the Stieffel invariants studied earlier and would seem to be
related to the ``moving frames'' of Cartan \cite{EC}.  Further, one
should be able to use the Lie algebra structure relative to the
Toeplitz Lie group to define characteristic forms which capture these
invariants.  Moreover, if one assumes that those invariants are the
same for the jet bundles for two line bundles, then the remaining
degrees of freedom in choosing the bundle map to be holomorphic
essentially amounts to a phase which in this case is a unitary-valued
function.  The existence question for such a phase would seem to be
related to the existence of complex structure and thus to Chern-Moser
invariants.

If one considers quotient modules for submodules of functions that vanish
to increasing order, then they form a natural inverse limit of Hilbert
modules whose limit will be $\mathcal M$.  In a dual manner, one
should be able to show that the direct limit of the jet bundles
constructed for these quotient modules has a limit equal to $\Omega
\times \mathcal M$.  One way of viewing these constructions would be
as expanding $\mathcal M$ as a ``Taylor series'' of modules over
$\scZ$.

Finally, assume that there is a global defining function $\phi$ for
$\scZ$ in $\scA(\Omega)$ and consider the operator it defines on the
quotient module defined by the functions which vanish to order $k$ in
the direction of $\scZ$.  Then $\phi$ defines a bundle map on the $k
\times k$ matrix-valued kernel Hilbert space for the quotient which
can be written as the scalar multiplier $\phi I_k$ plus a nilpotent
matrix-valued multiplier.  Such an operator can be seen to be
analogous to the spectral operators of Dunford \cite{DS3}.  That is
the case if one replaces a normal operator by a multiplication
operator on a space of holomorphic functions.  An abstract
characterization of operators having such a representation as well as
a study of their properties would seem to be of interest.

It is clear that the ideas and techniques of this paper
raise many questions that warrant additional study.

\begin{Acknowledgement} The authors thank V. Pati for
    many helpful discussions. For instance, the calculation of the
    second fundamental form in section \ref{pati} is entirely due to
    him. Also, the second named author would like to thank B. Bagchi
    for many useful conversations relating to the topic of this
    paper.
\end{Acknowledgement}

\bibliographystyle{amsplain}
\providecommand{\bysame}{\leavevmode\hbox to3em{\hrulefill}\thinspace}
\providecommand{\MR}{\relax\ifhmode\unskip\space\fi MR }
\providecommand{\MRhref}[2]{%
  \href{http://www.ams.org/mathscinet-getitem?mr=#1}{#2}
}
\providecommand{\href}[2]{#2}


\begin{thebibliography}{10}

\bibitem{Aron}
N.~Aronszajn, \emph{{Theory of reproducing kernels}}, Trans. Amer.Math. Soc.
  \textbf{68} (1950), 337--404.

\bibitem{BMJFA96}
B.~Bagchi and G.~Misra, \emph{{Homogeneous tuples of multiplication operators
  on twisted Bergman space}}, J. Funct. Anal. \textbf{136} (1996), 171 -- 213.

\bibitem{BMIAS01}
\bysame, \emph{{Homogeneous operators and projective representations of the
  M\"{o}bius group: a survey}}, Proc. Ind. Acad. Sc.(Math. Sci.) \textbf{111}
  (2001).

\bibitem{EC}
E.~Cartan, \emph{Lecons sur la g\'{e}om\'{e}trie des espaces de riemann}, Les
  Grands Classiques Gauthier-Villars, 1988, Reprint of the second (1946)
  edition.

\bibitem{chenrgd}
X.~Chen and R.~G. Douglas, \emph{Localization of {Hilbert} modules}, Mich.
  Math. J. \textbf{39} (1992), 443 -- 454.

\bibitem{CD}
M.~J. Cowen and R.~G. Douglas, \emph{Complex geometry and operator theory},
  Acta Math. \textbf{141} (1978), 187--261.

\bibitem{Bolyai}
\bysame, \emph{On operators possessing an open set of eigenvalues}, Memorial
  Conf. for F\'{e}jer-Riesz, Colloq. Math. Soc. J. Bolyai, 1980, pp.~323 --
  341.

\bibitem{CD85}
\bysame, \emph{Equivalence of {Connections}}, Adv. Math. \textbf{56} (1985), 39
  -- 91.

\bibitem{CS}
R.~E. Curto and N.~Salinas, \emph{Generalized {Bergman} kernels and the
  {Cowen-Douglas} theory}, Amer J. Math. \textbf{106} (1984), 447--488.

\bibitem{rgdgm}
R.~G. Douglas and G.~Misra, \emph{Geometric invariants for resolutions of
  {Hilbert} modules}, Operator Theory: Advances and Applications, Birkhauser,
  1993, pp.~83--112.

\bibitem{eqhm}
\bysame, \emph{Equivalence of quotient {Hilbert} modules}, Proc. Indian Acad.
  Sc. (Math. Sc.) \textbf{113} (2003), 281 -- 291.

\bibitem{qfie}
\bysame, \emph{Quasi-free resolutions of {Hilbert} modules}, Integr. Equ. Oper.
  Theory \textbf{99} (2003), 435 -- 456.

\bibitem{quasi}
\bysame, \emph{On quasi-free {Hilbert} modules}, New York J. Math., \textbf{11} (2005), 547 - 561. 

\bibitem{DMV}
R.~G. Douglas, G.~Misra, and C.~Varughese, \emph{On quotient modules - the case
  of arbitrary multiplicity}, J. Func. Anal. \textbf{174} (2000), 364--398.

\bibitem{rgdvip}
R.~G. Douglas and V.~I. Paulsen, \emph{{Hilbert modules over function
  algebras}}, Pitman research notes in mathematics, no. 217, Longman Scientific
  and Technical, 1989.

\bibitem{D-P-S-Y}
R.~G. Douglas, V.~I. Paulsen, C.-H. Sah, and K.~Yan, \emph{{Algebraic reduction
  and rigidity for Hilbert modules}}, Amer.\ J.\ Math. \textbf{117} (1995), 75
  -- 92.

\bibitem{DS3}
N.~Dunford and J.~T. Schwartz, \emph{{Linear operators. Part III. Spectral
  operators. With the assistance of William G. Bade and Robert G. Bartle}},
  John Wiley \& Sons, 1988, Reprint of the 1971 original.

\bibitem{FR}
S.~H. Ferguson and R.~Rochberg, \emph{{Higher-order Hilbert-Schmidt Hankel
  forms and tensors of analytic kernels}}, pre-print.

\bibitem{FR1}
\bysame, \emph{{Higher-order Hilbert-Schmidt Hankel forms and tensors of
  analytic kernels}}, Proc. London Math. Soc. \textbf{(3) 82} (2001), no.~1,
  110--130.

\bibitem{G-H}
P.~Griffiths and J.~Harris, \emph{{Principles of Algebraic Geometry}}, John
  Wiley \& Sons, 1978.

\bibitem{GR}
R.~C. Gunning and H.~Rossi, \emph{Analytic functions of several complex
  variables}, Prentice Hall, 1965.

  
\bibitem{AK} 
A.~Koranyi and G.~Misra, \emph{Homogeneous operators on 
Hilbert spaces of holomorphic functions  - I}, preprint.

\bibitem{KMisComptRendu} 
\bysame, {\em New construction of some  
homogeneous operators}, C. R. Acad. Sci. Paris, Ser. I 
\textbf{342} (2006), 933 -- 936. 

\bibitem{MSJOT90}
G.~Misra and N.~S.~N. Sastry, \emph{Homogeneous tuples of operators and
  holomorphic discrete series representation of some classical groups}, J.
  Operator Theory \textbf{24} (1990), 23 -- 32.

\bibitem{Peng-Zhang}
L.~Peng and G.~Zhang, \emph{Tensor products of holomorphic representations and
  bilinear differential operators}, J. Funct. Anal. \textbf{210} (2004), 171 --
  192.

\bibitem{Na-Fo}
B.~Sz-Nagy and C.~Foias, \emph{{Harmonic analysis of operators on Hilbert
  space}}, North Holland, 1970.

\bibitem{RW}
R.~O. Wells, \emph{Differential analysis on complex manifolds}, Springer, 1973.

\bibitem{wil}
D.~R. Wilkins, \emph{{Homogeneous vector bundles and Cowen-Douglas operators}},
  Intern. J. Math. \textbf{4} (1993), 503 -- 520.

\end{thebibliography}
\end{document}